# *Perelman's functional and reduced volume*

Hassan Jolany

*School of Mathematics, college of science, university of Tehran, Tehran, Iran*

hassan.jolany@khayam.ut.ac.ir

**Abstract:** In recent years, there has seen much interest and increased research activities on Perelman's paper. Section one and two of this paper aim to establish Perelman's local non-collapsing result for the Ricci flow. This will provide a positive lower bound on the injectivity radius for the Ricci flow under blow-up analysis. We also discuss the gradient flow formalism of the Ricci flow and Perelman's motivation from physics. In this paper, we survey some of the recent progress on Perelman's functional, reduced volume and reduced length.

# 1. Two functional $\mathcal{F}$ and $\mathcal{W}$ of Perelman

*In this section, we introduce two functional of Perelman $\mathcal{F}$ and $\mathcal{W}$, and discuss their relations with the Ricci flow. The functional $\mathcal{F}$ can be found in the literature on the String Theory, where it describes the flow energy effective action. It was not known whether the Ricci flow is a gradient flow until Perlman showed that the Ricci flow is, in a certain sense, the gradient flow of the functional $\mathcal{F}$. In the theory of dynamical systems we can take the two functional $\mathcal{F}$ and $\mathcal{W}$ as a Lyapunov type.*

### The $\mathcal{F}$ −functional

*Let $\mathcal{M}$ denote the space of smooth Riemannian metrics $g$ on $M$. We think of $\mathcal{M}$ formally as an infinite-dimensional manifold. Also $M$ is closed manifold and the tangent space $T_g\mathcal{M}$ consists of the symmetric covariant 2-tensors $v_{ij}$ on $M$. $f$ is a function $f\colon M \to \mathbb{R}$. Let $dv$ denote the Riemannian volume density associated to a metric $g$. The functional $\mathcal{F}\colon \mathcal{M} \times C^\infty(M) \to \mathbb{R}$ is given by*

$$\mathcal{F}(g,f) = \int_M (R + |\nabla f|^2)e^{-f}dV$$

Given $v_{ij} \in T_g M$ and $h = \delta f$. The evaluation of the differential $d\mathcal{F}$ on $(v_{ij}, h)$ is written as $\delta\mathcal{F}(v_{ij}, h)$. put $v = g^{ij}v_{ij}$.

*Theorem1.1 (Perelman) we have*

$$\delta\mathcal{F}(v_{ij}, h) = \int_M e^{-f}\left[-v_{ij}(R_{ij} + \nabla_i\nabla_j f) + \left(\frac{v}{2} - h\right)(2\Delta f - |\nabla f|^2 + R)\right]dV \quad (1.1)$$

Proof. First we prove $\delta R = -\Delta v + \nabla_i\nabla_j v_{ij} - R_{ij}v_{ij}$. In any normal coordinates at a fixed point, we have

$$\delta R_{ijl}^h = \frac{\partial}{\partial x^i}(\delta\Gamma_{jl}^h) - \frac{\partial}{\partial x^j}(\delta\Gamma_{il}^h)$$

$$= \frac{\partial}{\partial x^i}\left[\frac{1}{2}g^{hm}(\nabla_j v_{lm} + \nabla_l v_{jm} - \nabla_m v_{jl})\right] - \frac{\partial}{\partial x^j}\left[\frac{1}{2}g^{hm}(\nabla_i v_{lm} + \nabla_l v_{im} - \nabla_m v_{il})\right]$$

$$\delta R_{jl} = \frac{\partial}{\partial x^i}\left[\frac{1}{2}g^{im}(\nabla_j v_{lm} + \nabla_l v_{jm} - \nabla_m v_{jl})\right] - \frac{\partial}{\partial x^j}\left[\frac{1}{2}g^{im}(\nabla_i v_{lm} + \nabla_l v_{im} - \nabla_m v_{il})\right]$$

$$= \frac{1}{2}\frac{\partial}{\partial x^i}[\nabla_j v_l^i + \nabla_l v_j^i - \nabla^i v_{jl}] - \frac{1}{2}\frac{\partial}{\partial x^j}[\nabla_l v]$$

Therefore

$$\delta R = \delta(g^{il}R_{jl}) = -v_{jl}R_{jl} + g^{jl}\delta R_{jl} = -v_{jl}R_{jl} + \frac{1}{2}\frac{\partial}{\partial x^i}[\nabla^l v_l^i + \nabla_l v^{il} - \nabla^i v] - \frac{1}{2}\frac{\partial}{\partial x^j}[\nabla^j v]$$

$$= -v_{jl}R_{jl} + \nabla_i\nabla_l v_{il} - \Delta v$$

as $|\nabla f|^2 = g^{ij}\nabla_i f \nabla_j f$

We have $\delta|\nabla f|^2 = -v^{ij}\nabla_i f\nabla_j f + 2\langle\nabla f, \nabla h\rangle$, as $dV = \sqrt{\det(g)}\, dx_1 \ldots dx_n$, we have $\delta(dV) = \frac{v}{2}dV$, so

$$\delta(e^{-f}dV) = \left(\frac{v}{2} - h\right)e^{-f}dV \quad (1.2)$$

Putting this together gives

$$\delta\mathcal{F} = \int_M e^{-f}\left[-\Delta v + \nabla_i\nabla_j v_{ij} - R_{ij}v_{ij} - v_{ij}\nabla_i f\nabla_j f + 2\langle\nabla f, \nabla h\rangle + (R + |\nabla f|^2)\left(\frac{v}{2} - h\right)\right]dV$$

(1.3)

The goal now is to rewrite the right-hand side of (1.3) so that $v_{ij}$ and $h$ appear algebraically, i.e. without derivatives. as

$$\Delta e^{-f} = (|\nabla f|^2 - \Delta f)e^{-f}$$

we have

$$\int_M e^{-f}[-\Delta v]dV = -\int_M (\Delta e^{-f})v dV = \int_M e^{-f}(\Delta f - |\nabla f|^2)v dV$$

next

$$\int_M e^{-f}\nabla_i\nabla_j v_{ij} dV = \int_M (\nabla_i\nabla_j e^{-f})v_{ij}dV$$

$$= -\int_M \nabla_i(e^{-f}\nabla_j f)v_{ij}dV = \int_M e^{-f}(\nabla_i f\nabla_j f - \nabla_i\nabla_j f)v_{ij}dV$$

finally

$2 \int_M e^{-f} \langle \nabla f, \nabla h \rangle \, dV = -2 \int_M \langle \nabla e^{-f}, \nabla h \rangle \, dV = 2 \int_M (\Delta e^{-f}) h \, dV = 2 \int_M e^{-f} (|\nabla f|^2 - \Delta f) h \, dV$

then

$$\delta \mathcal{F} = \int_M e^{-f} \left[ \left(\frac{v}{2} - h\right)(2\Delta f - 2|\nabla f|^2) - v_{ij}(R_{ij} + \nabla_i \nabla_j f) + \left(\frac{v}{2} - h\right)(R + |\nabla f|^2) \right] dV$$

$$= \int_M e^{-f} \left[ -v_{ij}(R_{ij} + \nabla_i \nabla_j f) + \left(\frac{v}{2} - h\right)(2\Delta f - |\nabla f|^2 + R) \right] dV$$

So proof is complete ∎

*Remark 1.1(see [4-5]) notice that*

1) the functional $\mathcal{F}$ is invariant under diffeomorphism i.e. $\mathcal{F}(M, \phi^* g, f \circ \phi) = \mathcal{F}(M, g, \phi)$ for any diffeomorphisms $\phi$
2) Also for any $c > 0$ and $b$, $\mathcal{F}(M, c^2 g, f + b) = c^{n-2} e^{-b} \mathcal{F}(M, g, f)$.

Example: Let $(M, g)$ be Euclidean space $M = R^n$ and let

$$f(t, x) = \frac{|x|^2}{4T} + \frac{n}{2} \log 4\pi T = -\log \left[ (4\pi T)^{\frac{-n}{2}} e^{\frac{-|x|^2}{4T}} \right]$$

Where $\tau = t_0 - t$, notice that $e^{-f} dx$ is the Gaussian measure, which solves the backward heat equation. If $t < t_0$. This choice of $g$ and $f$ satisfy the equations

$$\frac{\partial g}{\partial t} = -2Rc(g)$$
$$\frac{\partial f}{\partial t} = -\Delta f - R + |\nabla f|^2$$

we can check

$$\frac{\partial f}{\partial t} = \frac{\partial}{\partial t} \left[ \frac{|x|^2}{4T} + \frac{n}{2} \log 4\pi T \right] = \frac{|x|^2}{4T^2} - \frac{n}{2T}$$

and

$$\nabla f = \frac{x}{2T}$$

so $|\nabla f|^2 = \frac{|x|^2}{4T^2}$ and $\Delta f = \frac{n}{2T}$.

also because we know $\int_{R^n} e^{\frac{-|x|^2}{4T}} dV = (4\pi T)^{\frac{n}{2}}$. By differentiating with respect to $\tau$ gives

$$\int_{R^n} \frac{|x|^2}{4T^2} e^{\frac{-|x|^2}{4T}} dV = (4\pi T)^{\frac{n}{2}} \frac{n}{2T}.$$

so

$$\int_{R^n} |\nabla f|^2 e^{-f} dv = \frac{n}{2T}$$

Then $\mathcal{F}(t) = \frac{n}{2T} = \frac{n}{2(t_0 - t)}$. In particular, this is non-decreasing as a function of $t \in [0, t_0)$ ∎

Theorem 1.2 let $g_{ij}$ and $f(t)$ evolve according to the coupled flow

$$\frac{\partial g_{ij}}{\partial t} = -2R_{ij}$$
$$\frac{\partial f}{\partial t} = -\Delta f - R + |\nabla f|^2$$

Then $\int_M e^{-f} dV$ is constant.

Proof. By applying the chain rule we have

$$\frac{\partial}{\partial t} dV = \frac{\partial}{\partial t}\sqrt{det g_{ij}} dx^1 \wedge dx^2 \wedge ... \wedge dx^n$$

$$= \frac{1}{2\sqrt{det g_{ij}}} g^{ij}(-2R_{ij}) det g dx^1 \wedge dx^2 \wedge ... \wedge dx^n = -RdV$$

because we have $\frac{d}{dt} det A = (A^{-1})^{ij}\left(\frac{dA_{ij}}{dt}\right) det A$

hence

$$\frac{d}{dt}(e^{-f} dV) = e^{-f}\left(-\frac{\partial f}{\partial t} - R\right) dV = (\Delta f - |\nabla f|^2) e^{-f} dV = -\Delta(e^{-f}) dV$$

because $\Delta(e^{-f}) = (|\nabla f|^2 - \Delta f)e^{-f}$.

So it then follows that

$$\frac{d}{dt}\int_M e^{-f} dV = -\int_M \Delta(e^{-f}) dV$$

Because M is closed manifold according to divergence theorem

$$\int_M \Delta(e^{-f}) dV = 0$$

So this finishes the proof of the theorem ∎

we would like to get rid of the $\left(\frac{v}{2} - h\right)(2\Delta f - |\nabla f|^2 + R)$ term in (1.1). We can do this by restricting our variations so that $\frac{v}{2} - h = 0$. From (1.2), this amounts to assuming that assuming $e^{-f} dV$ is fixed. We now fix a smooth measure $dm$ on $M$ and relate $f$ to $g$ by requiring that $e^{-f} dV = dm$. Equivalently, we define a section $s : \mathcal{M} \to \mathcal{M} \times C^\infty(M)$ by $s(g) = \left(g, \ln\left(\frac{dV}{dm}\right)\right)$. Then the composition $\mathcal{F}^m = \mathcal{F} \circ s$ is a function on $\mathcal{M}$ and its deferential is given by

$$d\mathcal{F}^m(v_{ij}) = \int_M e^{-f}[-v_{ij}(R_{ij} + \nabla_i \nabla_j f)] dV \quad (1.4)$$

Defining a formal Riemannian metric on $\mathcal{M}$ by

$$\langle v_{ij}, v_{ij}\rangle_g = \frac{1}{2}\int_M v^{ij} v_{ij} dm$$

The gradient flow of $\mathcal{F}^m$ on $\mathcal{M}$ is given by

$$\frac{\partial}{\partial t} g_{ij} = -2(R_{ij} + \nabla_i \nabla_j f) \quad (1.5)$$

The induced flow equation for $f$ is

$$\frac{\partial f}{\partial t} = \frac{1}{2} g^{ij} \frac{\partial}{\partial t} g_{ij} = -\Delta f - R$$

As with any gradient flow, the function $\mathcal{F}^m$ is non-decreasing along the flow line with its derivative being given by the length squared of the gradient, i.e.

$$\frac{\partial}{\partial t}\mathcal{F}^m = 2\int_M |R_{ij} + \nabla_i \nabla_j f|^2 dm$$

We can check that if $g_{ij}(t)$ and $f(t)$ evolve according to the coupled flow

$$\begin{cases} \frac{\partial g_{ij}}{\partial t} = -2R_{ij} \\ \frac{\partial f}{\partial t} = -\Delta f + |\nabla f|^2 - R \end{cases} \quad (1.7)$$

Then $\frac{d}{dt}\mathcal{F}(g_{ij}(t), f(t)) = 2\int_M |R_{ij} + \nabla_i\nabla_j f|^2 e^{-f}\, dV$. Now we prove that to obtaining a solution of (1.5) and (1.6) is to show that it is some how equivalent to the decoupled system of (1.7) ∎

Proposition1.1 (see [5]) defining $\hat{g}(t) = \sigma(t)\psi_t^*(g(t))$, we have

$$\frac{\partial \hat{g}}{\partial t} = \sigma'(t)\psi_t^*(g) + \sigma(t)\psi_t^*\left(\frac{\partial g}{\partial t}\right) + \sigma(t)\psi_t^*(\mathcal{L}_X g) \quad (1.9)$$

Also for some function $f: M \to R$, we then have

$$\mathcal{L}_{(\nabla f)} g = 2\text{Hess}(f)$$

Theorem1.3 the solutions of (1.5) and (1.6) may be generated by pulling back solutions of (1.7) by an appropriate time-dependent diffeomorphism.

Proof. By previous proposition we know $\mathcal{L}_X g = -2\text{Hess}(f)$ where $X(t) = -\nabla f$. We fix $\sigma(t) \equiv 1$. Now we define $\hat{g}(t)$ by (1.8), it will evolve, by (1.9), according to

$$\frac{\partial \hat{g}}{\partial t} = \psi_t^*(t)[-2\text{Ric}(g) - 2\text{Hess}_g(f)]$$

Where $\text{Hess}_g(f)$ is the Hessian of $f$ with respect to the metric $g$. Keeping in mind that $\psi_t: (M, \hat{g}) \to (M, g)$ is an isometry, we may then write $\hat{f} := f \circ \psi_t$ to give

$$\frac{\partial \hat{g}}{\partial t} = -2(\text{Ric}(\hat{g}) + \text{Hess}_{\hat{g}}(\hat{f}))$$

The evolution of $\hat{f}$ is then found by the chain rule

$$\frac{\partial \hat{f}}{\partial t}(x, t) = \frac{\partial f}{\partial t}(\psi_t(x), t) + X(f)(\psi_t(x), t)$$
$$= [(-\Delta_g f + |\nabla f|^2 - R_g) - |\nabla f|^2](\psi_t(x), t) = [-\Delta_{\hat{g}} \hat{f} - R_{\hat{g}}](x, t)$$

(Here we have used of this fact that $\mathcal{L}_{\nabla f} f = |\nabla f|^2$)

Thus we have a solution to $\frac{\partial \hat{f}}{\partial t} = -\Delta_{\hat{g}} \hat{f} - R_{\hat{g}}$ ∎

We would like to use this to develop a controlled quantity for Ricci flow, but we need to eliminate $f$. This can be accomplished by taking an infimum, defining

$$\lambda(M, g) = \inf_{f: \int_M e^{-f} dV = 1} \mathcal{F}(M, g, f)$$

Lemma1.1 (see [5]) let $g_{ij}(t)$ and $f(t)$ evolve according to the coupled flow

$$\begin{cases} \frac{\partial g_{ij}}{\partial t} = -2R_{ij} \\ \frac{\partial f}{\partial t} = -\Delta f + |\nabla f|^2 - R \end{cases}$$

Then $\mathcal{F}(g_{ij}(t), f(t))$ is non-decreasing in time and monotonicity is strict unless we are on a steady gradient solution.

Proof. By applying previous computations we can show that

$$\frac{\partial}{\partial t}\mathcal{F}(g_{ij}(t),f(t)) = 2\int_M |R_{ij}+\nabla_i\nabla_j f|^2 e^{-f}\,dV$$

So proof is complete ∎

So by previous lemma we obtain

$$\lambda(g_{ij}(t)) \leq \mathcal{F}(g_{ij}(t),f(t)) \leq \mathcal{F}(g_{ij}(t_0),f(t_0)) = \lambda(g_{ij}(t_0)) \quad (1.10)$$

For $t < t_0$ and $\int_M e^{-f}dV = 1$.

Definition1.1 a steady breather is a Ricci flow solution on an interval $[t_1, t_2]$ that satisfies the equation $g(t_2) = \phi^* g(t_1)$ for some $\phi \in Diff(M)$.

Now we show that a steady breather on a compact manifold is a gradient steady soliton.

Theorem1.4 (Perelman[5]) a steady breather is a gradient steady soliton.

Proof. We have $\lambda(g(t_2)) = \lambda(\phi^* g(t_1)) = \lambda(g(t_1))$. Because $\lambda$ is invariant under diffeomorphism. So by (1.10), $\mathcal{F}(g(t),f(t))$ must be constant in t. From

$$\frac{\partial}{\partial t}\mathcal{F}(g_{ij}(t),f(t)) = 2\int_M |R_{ij}+\nabla_i\nabla_j f|^2 e^{-f}\,dV$$

We conclude, $R_{ij}+\nabla_i\nabla_j f = 0$. Then $R + \Delta f = 0$ and so the system

$$\begin{cases} \frac{\partial g_{ij}}{\partial t} = -2R_{ij} \\ \frac{\partial f}{\partial t} = -\Delta f + |\nabla f|^2 - R \end{cases}$$

Change to

$$\frac{\partial g_{ij}}{\partial t} = -2R_{ij}$$
$$\frac{\partial f}{\partial t} = |\nabla f|^2$$

This is a gradient expanding soliton ∎

Lemma 1.2(see [4-5]) we have

$$\frac{\partial \lambda}{\partial t} \geq \frac{2}{n}\lambda^2(t)$$

Definition1.2 an expanding breather is a Ricci flow solution on an interval $[t_1, t_2]$ that satisfies the equation $g(t_2) = c\phi^* g(t_1)$ for some $c > 1$ and $\phi \in Diff(M)$.

Proposition1.2 an expanding breather is a gradient expanding solution.

Lemma1.3 (see [8]) $\lambda(M, g)$ is finite.

Lemma 1.4(see [8]) $\lambda(M, g)$ is the least number for which one has the inequality

$$\int_M 4|\nabla u|_g^2 + R|u|^2 dV \geq \lambda(M,g)\int_M |u|^2 dV$$

For all $u$ in the sobolev space $H^1(M)$. (Note $\|f\|_{H^1} = \int (|\nabla f|_g^2 + f^2)dV$ for $C^1$ functions)

### *The $\mathcal{W}$ − Functional*

*We know that the metric $g_{ij}(t)$ evolving by Ricci flow is called a breather, if for some $t_1 < t_2$ and $\alpha > 0$ the metrics $\alpha g_{ij}(t_1)$ and $g_{ij}(t_2)$ differ only by 0 diffeomorphism; the cases $\alpha = 1, \alpha < 1, \alpha > 1$ correspond to steady, shrinking and expanding breathers, respectively. In order to handle the shrinking case when $\lambda(M,g) > 0$ we need to replace our functional $\mathcal{F}$ by its generalization, which contains explicit insertions of the scale parameter, to be denoted by $\tau$. Thus consider the functional.*

$$\mathcal{W}(g_{ij}, f, \tau) = \int_M [\tau(|\nabla f|^2 + R) + f - n](4\pi\tau)^{\frac{-n}{2}} e^{-f} dV$$

*Restricted to $f$ satisfying*

$$\int_M (4\pi\tau)^{\frac{-n}{2}} e^{-f} dV = 1$$

*$\tau > 0$. Where $g_{ij}$ is a Riemannian metric, $f$ is a smooth function on $M$ and $\tau$ is a positive scale parameter. Also for any positive number $a$ and any diffeomorphism $\varphi$*

$$\mathcal{W}(a\varphi^* g_{ij}, \varphi^* f, a\tau) = \mathcal{W}(g_{ij}, f, \tau).$$

*Now we start with first variation for $\mathcal{W}$.*

*Theorem 1.5 (Perelman) assume that $\delta g_{ij} = v_{ij}$ and $\delta f = h$. Put $\sigma = \delta\tau$ then we have*

$$\delta\mathcal{W}(v_{ij}, h, \sigma) = \int_M \left[\sigma(R + |\nabla f|^2) - \tau v_{ij}(R_{ij} + \nabla_i \nabla_j f) + h + [\tau(2\Delta f - |\nabla f|^2 + R) + f - n]\left(\frac{v}{2} - h - \frac{n\sigma}{2\tau}\right)\right](4\pi\tau)^{\frac{-n}{2}} e^{-f} dV$$

*Proof. By $\frac{\partial}{\partial t} dV = \frac{v}{2} dV$. We see*

$$\delta\left((4\pi\tau)^{\frac{-n}{2}} e^{-f} dV\right) = \left(\frac{v}{2} - h - \frac{n\sigma}{2\tau}\right)(4\pi\tau)^{\frac{-n}{2}} e^{-f} dV$$

*but we know*

$$\delta\mathcal{F}(v_{ij}, h) = \int_M e^{-f} \left[-v_{ij}(R_{ij} + \nabla_i \nabla_j f) + \left(\frac{v}{2} - h\right)(2\Delta f - |\nabla f|^2 + R)\right] dV$$

*we obtain*

$$\delta\mathcal{W} = \int_M \left[\sigma(R + |\nabla f|^2) + \tau\left(\frac{v}{2} - h\right)(2\Delta f - 2|\nabla f|^2) - \tau v_{ij}(R_{ij} + \nabla_i \nabla_j f) + h \right.$$
$$\left. + [\tau(R + |\nabla f|^2) + f - n]\left(\frac{v}{2} - h - \frac{n\sigma}{2\tau}\right)\right](4\pi\tau)^{\frac{-n}{2}} e^{-f} dV$$

*also we have*

$$\Delta e^{-f} = (|\nabla f|^2 - \Delta f)e^{-f}$$

*therefore*

$$\delta\mathcal{W}(v_{ij}, h, \sigma) = \int_M \left[\sigma(R + |\nabla f|^2) - \tau v_{ij}(R_{ij} + \nabla_i \nabla_j f) + h \right.$$
$$\left. + [\tau(2\Delta f - |\nabla f|^2 + R) + f - n]\left(\frac{v}{2} - h - \frac{n\sigma}{2\tau}\right)\right](4\pi\tau)^{\frac{-n}{2}} e^{-f} dV$$

*So proof is complete* ■

*Definition1.3 the arguments $g, f$ and $\tau$ are called compatible if*
$$\int_M \frac{e^{-f}}{(4\pi\tau)^{\frac{n}{2}}} dV = 1$$

*Lemma1.5(see[5]) under the transformation $(g, f, \tau) \to (\tau^{-1}g, f, 1)$ compatibility is preserved, as is the functional :*
$$\mathcal{W}(g, f, \tau) = \mathcal{W}(\tau^{-1}g, f, 1)$$

*Note that on $R^n$, the Gaussian measure $d\mu$ is defined in terms of the lebesgue measure $dx$ by*
$$d\mu = (2\pi)^{\frac{-n}{2}} e^{\frac{-|x|^2}{2}} dx$$

*The normalization being chosen so that*
$$\int_{R^n} d\mu = 1$$

*Lemma 1.6(see [7])(L.Gross) If $v: R^n \to R$ is, say, smooth and satisfies $v, |\nabla v| \in L^2(d\mu)$ then*
$$\int v^2 \ln|v| d\mu \leq \int |\nabla v|^2 d\mu + \left(\int v^2 d\mu\right) \ln \left(\int v^2 d\mu\right)^{\frac{1}{2}}$$

*so if we choose $v$ so that $\int v^2 d\mu = 1$, then inequality becomes*
$$\int v^2 \ln|v| d\mu \leq \int |\nabla v|^2 d\mu$$

*That called log-sobolev inequality.*

*Theorem1.6(see[5]) (Perelman) Let $g$ denote the flat metric on $R^n$. if $f$ and $\tau$ are compatible with $g$, then*
$$\mathcal{W}(g, f, \tau) \geq 0$$

*Proof let $f: R^n \to R$ be compatible with $g$ and $\tau$, which in this situation means that*
$$\int_{R^n} \frac{e^{-f}}{(4\pi\tau)^{\frac{n}{2}}} dx = 1$$

*If we set $\tau = \frac{1}{2}$ we obtain $\int_{R^n} \frac{e^{-f}}{(2\pi)^{\frac{n}{2}}} dx = 1$. If we define $v = e^{\frac{|x|^2}{4} - \frac{f}{2}}$, we have*
$$v^2 d\mu = e^{\frac{|x|^2}{2} - f} (2\pi)^{\frac{-n}{2}} e^{\frac{-|x|^2}{2}} dx = (2\pi)^{\frac{-n}{2}} e^{-f} dx$$

*so $\int_{R^n} v^2 dx = 1$. Therefore, by the log-sobolev inequality we obtain, $\int v^2 \ln|v| d\mu \leq \int |\nabla v|^2 d\mu$.*

*by computing the left-hand side and right hand side of this inequality we get*
$$\int v^2 \ln|v| d\mu = \int e^{\frac{|x|^2}{2} - f} \left(\frac{|x|^2}{4} - \frac{f}{2}\right) (2\pi)^{\frac{-n}{2}} e^{\frac{-|x|^2}{2}} dx$$
$$= \int \left(\frac{|x|^2}{4} - \frac{f}{2}\right) \frac{e^{-f}}{(2\pi)^{\frac{n}{2}}} dx$$

also $\nabla v = \left(\frac{x}{2} - \frac{\nabla f}{2}\right) e^{\frac{|x|^2}{4} - \frac{f}{2}}$

which gives us

$$|\nabla v|^2 = \left(\frac{|x|^2}{4} - \frac{x.\nabla f}{2} + \frac{|\nabla f|^2}{4}\right) e^{\frac{|x|^2}{2} - f}$$

therefore

$$|\nabla v|^2 d\mu = \left(\frac{|x|^2}{4} - \frac{x.\nabla f}{2} + \frac{|\nabla f|^2}{4}\right) \frac{e^{-f}}{(2\pi)^{\frac{n}{2}}} dx$$

the $integration - by - parts$ formula gives us

$$-\int \frac{x.\nabla f}{2} \frac{e^{-f}}{(2\pi)^{\frac{n}{2}}} dx = \frac{1}{2}\int x.\nabla(e^{-f}) \frac{dx}{(2\pi)^{\frac{n}{2}}}$$

but we can compute $\nabla.x = n$, so

$$-\int \frac{x.\nabla f}{2} \frac{e^{-f}}{(2\pi)^{\frac{n}{2}}} dx = \frac{-n}{2}\int \frac{e^{-f}}{(2\pi)^{\frac{n}{2}}} dx$$

so

$$\int |\nabla v|^2 d\mu = \int \left(\frac{|x|^2}{4} - \frac{n}{2} + \frac{|\nabla f|^2}{4}\right) \frac{e^{-f}}{(2\pi)^{\frac{n}{2}}} dx$$

and the log-sobolev inequality gives us

$$\int \left(\frac{|x|^2}{4} - \frac{f}{2}\right) \frac{e^{-f}}{(2\pi)^{\frac{n}{2}}} dx \leq \int \left(\frac{|x|^2}{4} - \frac{n}{2} + \frac{|\nabla f|^2}{4}\right) \frac{e^{-f}}{(2\pi)^{\frac{n}{2}}} dx$$

so

$$\mathcal{W}\left(g, f, \frac{1}{2}\right) = \int \left[\frac{1}{2}|\nabla f|^2 + f - n\right] \frac{e^{-f}}{(2\pi)^{\frac{n}{2}}} dx \geq 0$$

by the scale invariance $\mathcal{W}(g, f, \tau) = \mathcal{W}\left(\frac{1}{2\tau}g, f, \frac{1}{2}\right)$ and because $(R^n, g)$ is preserved under the homothetic scaling, we conclude

$$\mathcal{W}(g, f, \tau) \geq 0$$

and proof will be complete■

Remark 1.2(see [4]) to easily we can check that for any $f$, $\tau$ compatible with $g$, $\mathcal{W}(g, f, \tau) = 0$ if and only if $f(x) \equiv \frac{x^2}{4\tau}$.

We have a analogous theorem for $\mathcal{W}$, like $\mathcal{F}$. We see that $\mathcal{W}$, is increasing under the Ricci flow when $f$ and $\tau$ are made to evolve appropriately.

Theorem1.7 (Perelman) If $g_{ij}(t), f(t)$ and $\tau(t)$ evolve according to the system

$$\begin{cases} \frac{\partial g_{ij}}{\partial t} = -2R_{ij} \\ \frac{\partial f}{\partial t} = -\Delta f + |\nabla f|^2 - R + \frac{n}{2\tau} \\ \frac{\partial \tau}{\partial t} = -1 \end{cases}$$

*Then we have the identity*

$$\frac{d}{dt}\mathcal{W}(g_{ij}(t), f(t), \tau(t)) = \int_M 2\tau \left| R_{ij} + \nabla_i\nabla_j f - \frac{1}{2\tau}g_{ij} \right|^2 (4\pi\tau)^{-\frac{n}{2}} e^{-f} dV$$

*and $\int_M (4\pi\tau)^{-\frac{n}{2}} e^{-f} dV$ is constant .In Particular .in particular $\mathcal{W}(g_{ij}(t), f(t), \tau(t))$ is non-decreasing in time and monotonicity is strict unless we are on a shrinking gradient soliton.*

*Now we give an example of a gradient shrinking soliton .*

*Example consider $R^n$ with the flat metric, constant in time $t \in (-\infty, 0)$ and let $\tau = -t$ and $f(\tau, x) = \frac{x^2}{4\tau}$.*

*proof. Because $f(\tau, x) = \frac{x^2}{4\tau}$ so we get $e^{-f} = e^{-\frac{x^2}{4\tau}}$. To easily we can check $(g(t), f(t), \tau(t))$ satisfies the following system*

$$\begin{cases} \frac{\partial g_{ij}}{\partial t} = -2R_{ij} \\ \frac{\partial f}{\partial t} = -\Delta f + |\nabla f|^2 - R + \frac{n}{2\tau} \\ \frac{\partial \tau}{\partial t} = -1 \end{cases}$$

*and $\int (4\pi\tau)^{-\frac{n}{2}} e^{-f} dV = 1$*
*now $\tau(|\nabla f|^2 + R) + f - n = \tau \frac{|x|^2}{4\tau^2} + \frac{|x|^2}{4\tau} - n = \frac{|x|^2}{2\tau} - n$*
*so it follows from of this fact that we have*

$$\int_{R^n} e^{-\frac{|x|^2}{4\tau}} dV = (4\pi\tau)^{\frac{n}{2}}$$

*and*

$$\int e^{-\frac{|x|^2}{4\tau}} \frac{|x|^2}{4\tau^2} dV = (4\pi\tau)^{\frac{n}{2}} \frac{n}{2\tau}$$

*Therefore $\mathcal{W}(t) = 0$ for all t.So proof is complete* ∎

*Remark1.3 Now if $g$ is the Euclidean metric and we let $u = (4\pi\tau)^{-\frac{n}{2}} e^{-f}$, we see that*

$$\log u = -\frac{n}{2}\log(4\pi\tau) - f$$
$$|\nabla u|^2 = (4\pi\tau)^{-n}|\nabla f|^2 e^{-2f}$$

*so* $\quad |\nabla f|^2 = \frac{|\nabla u|^2}{u^2}$

*Therefore*

$$\mathcal{W}(M, g, f, \tau) = \int \left[\tau \frac{|\nabla u|^2}{u^2} - u\log u\right] dx - \frac{n}{2}\log(4\pi\tau) - n$$

*but we know that $\mathcal{W} \geq 0$ so it implies a log-sobolev inequality*

$$\tau \int \frac{|\nabla u|^2}{u^2} dx \geq \int u\log u \, dx + \frac{n}{2}(4\pi\tau) + n$$

*If we set $\phi^2 = u$*

$$4\tau \int |\nabla\phi|^2 dx \geq \frac{1}{\tau}\int \phi^2 \log\phi^2 dx + \frac{n}{2}\log(4\pi\tau) + n$$

For the general case, we have

$$\mathcal{W}(M,g,f,\tau) = \int \left[\tau\left(Ru + \frac{|\nabla u|^2}{u^2}\right) - u\log u\right] dV - \frac{n}{2}\log(4\pi\tau) - n \quad (1.11)$$

one can show that

$$\mathcal{W}(M,g,f,\tau) \geq -c(M,g,\tau)$$

so according to (1.11) we get

$$\tau \int R\phi^2 dV + \tau \int 4|\nabla\phi|^2 dV \geq -c + \int \phi^2 \log\phi^2 dV + \frac{n}{2}\log(4\pi\tau) + n \quad (1.12)$$

*Definition 1.4* We set

$$\mu(g_{ij},\tau) = \inf\left\{\mathcal{W}(g_{ij},f,\tau): f \in C^\infty(M), \frac{1}{(4\pi\tau)^{\frac{n}{2}}}\int_M e^{-f} dV = 1\right\}$$

Which according to (1.12) is the best possible constant $-c$.

*Remark 1.4 (see [4])* $\mu$ is finite.

*Definition 1.5* a shrinking breather is a Ricci flow solution on $[t_1, t_2]$ that satisfies $g(t_2) = c\phi^* g(t_1)$ for some $c < 1$ and $\phi \in Diff(M)$.

*Definition 1.6* a shrinking soliton lives on a time interval $(-\infty, 0)$. a gradient shrinking soliton satisfies the equations

$$\frac{\partial g_{ij}}{\partial t} = -2R_{ij} = 2\nabla_i \nabla_j f + \frac{g_{ij}}{t}$$

$$\frac{\partial f}{\partial t} = |\nabla f|^2$$

*Remark 1.5 (see [5])* one can show that for any time $t \leq t_0$ we have

$$\mu(M,g(t),\tau(t)) \leq \mathcal{W}(M,g(t),f(t),\tau(t)) \leq \mu(M,g(t_0),\tau(t_0)) \quad (1.13)$$

*Theorem 1.8* A shrinking breather is a gradient shrinking soliton.

*Proof* put $t_0 = \frac{t_2 - ct_1}{1-c}$. Then if $\tau_1 = t_0 - t_1$ and $\tau_2 = t_0 - t_2$, we get $\tau_2 = c\tau_1$. Therefore because the functional $\mathcal{W}$ is invariant under simultaneous scaling of $\tau$ and $g_{ij}$ are invariant under diffeomorphism, so

$$\mu(g(t_2),\tau_2) = \mu\left(\frac{\tau_2}{\tau_1}\phi^* g(t_1), \tau_2\right) = \mu(\phi^* g(t_1), \tau_1) = \mu(g(t_1), \tau_1)$$

so by (1.13) and this fact that

$$\frac{d}{dt}\mathcal{W}(g_{ij}(t), f(t), \tau(t)) = \int_M 2\tau \left|R_{ij} + \nabla_i \nabla_j f - \frac{1}{2\tau}g_{ij}\right|^2 (4\pi\tau)^{\frac{-n}{2}} e^{-f} dV$$

it follows that the solution is a gradient shrinking soliton.

*Remark 1.6 (see [8]):* $\mu(g_{ij}(t), \tau - t)$ is non-decreasing along the Ricci flow.

*Proposition 1.3 (see [8])* $\mu(g,\tau)$ is negative for small $\tau > 0$ and tends to zero as $\tau \to 0$.

# Recent developments on Perelman's functional $\mathcal{F}$ and $\mathcal{W}$

*Definition1.7* In [2] Jun-Fang Li introduced the following $\mathcal{F}$ −functional,

$$\mathcal{F}_k(g,f) = \int_M (kR + |\nabla f|^2)e^{-f} d\mu$$

Where $k \geq 1$. when $k = 1$, this is the $\mathcal{F}$ −functional.

*The following theorem is analogous result, like, $\mathcal{F}$.*

*Theorem1.9 (see [2]) suppose the Ricci flow of $g(t)$ exists for $[0,T)$, then all the functional $\mathcal{F}_k(g,f)$ will be monotone under the following coupled system, i.e.*

$$\begin{cases} \frac{\partial g_{ij}}{\partial t} = -2R_{ij} \\ \frac{\partial f}{\partial t} = -\Delta f + |\nabla f|^2 - R \end{cases}$$

$$\frac{d}{dt}\mathcal{F}_k(g_{ij},f) = 2(k-1)\int_M |Rc|^2 e^{-f} d\mu + 2\int_M |R_{ij} + \nabla_i \nabla_j f|^2 e^{-f} d\mu \geq 0$$

define $\lambda_k(g) = \inf \mathcal{F}_k(g,f)$, where infimum is taken over all smooth $f$, satisfying $\int_M e^{-f} d\mu = 1$. and we assume $\lambda_1(g) = \lambda(g)$.

*Theorem1.10(see [2]and[11]) $\lambda(g)$ is the lowest eigenvalue of the parameter $-4\Delta + R$ and the non-decreasing of the $\mathcal{F}$ functional implies the non-decreasing of $\lambda(g)$. as an application, Perelman was able to show that there is no non-trivial steady or expanding Ricci breathers on closed manifolds.*

*Theorem1.11(see[11]) on a compact Riemannian manifold $(M, g(t))$, where $g(t)$ satisfies the Ricci flow equation for $t \in [0, T)$, the lowest eigenvalue $\lambda_k$ of the operator $-4\Delta + kR$ is non-decreasing under the Ricci flow. The monotonicity is strict unless the metric is Ricci-flat.*

*X. D. Cao considered the eigenvalues of the operator*

$$-\Delta + \frac{R}{2}$$

*On manifolds with nonnegative cuvvature operator. He showed that the eigenvalues of these manifolds are non-decreasing along the Ricci flow.*

*Corollary1.1 On a compact Riemannian manifold, the lowest eigenvalues of the operator $-\Delta + \frac{R}{2}$ are non-decreasing under the Ricci flow .*

*proof let $k = 2$, then $\frac{1}{4}\lambda_2$ is the lowest eigenvalue of $-\Delta + \frac{R}{2}$ and the result will follows* ■

*Theorem1.12 (see [11]) Let $g(t), t \in [0,T)$, be a solution to the Ricci flow on a closed Riemannian manifold $M^n$. Assume that there is a $C^1$-family of smooth function $f(t) > 0$, which satisfy*

$$\lambda(t)f(t) = -\Delta_{g(t)}f(t) + \frac{1}{2}R_{g(t)}f(t)$$

$$\int f^2(t)d\mu_{g(t)} = 1$$

Where $\lambda(t)$ is a function of $t$ only. Then

$$2\frac{d}{dt}\lambda(t) = 4\int R_{ij}\nabla^i f \nabla^j f\, d\mu + 2\int |Rc|^2 f^2\, d\mu$$

$$= \int |R_{ij} + \nabla_i\nabla_j\varphi|^2 e^{-\varphi}d\mu + \int |Rc|^2 e^{-\varphi}\, d\mu \geq 0$$

## Entropy functional for diffusion operator

Let $(M,g)$ be a compact Riemannian manifold, $\phi \in C^2(M)$. Let
$$L = \Delta - \nabla\phi.\nabla, \quad d\mu = e^{-\phi}dV$$

Let
$$u = \frac{e^{-f}}{(4\pi t)^{\frac{m}{2}}}$$

be a positive solution of
$$(\partial_t - L)u = 0$$

Inspired by the work of Perelman and Ni, we have the following results.

Theorem 1.13(X.-D. Li 2006) let

$$H_m(u,t) = \int_M u\log u\, d\mu - \left(\frac{m}{2}\log(4\pi t) + \frac{m}{2}\right)$$

$$\mathcal{W}(u,t) = \int_M (t|\nabla f|^2 + f - m)\frac{e^{-f}}{(4\pi t)^{\frac{m}{2}}}d\mu$$

Then

$$\frac{d}{dt}H_m(u,t) = -\int_M \left(L\log u + \frac{m}{2t}\right)u\, d\mu$$

$$\mathcal{W}(u,t) = \frac{d}{dt}(tH_m(u,t))$$

Theorem 1.14 (X.-D. Li) let u be a positive solution of the heat equation
$$\left(\frac{\partial}{\partial t} - L\right)u = 0$$

Suppose that

$$Ric_{m,n}(L) := Ric + \nabla^2\phi - \frac{\nabla\phi \otimes \nabla\phi}{m-n} \geq 0$$

Then

$$\frac{d\mathcal{W}(u,t)}{dt} = -2\int_M \tau\left(\left|\nabla^2 f - \frac{g}{2\tau}\right|^2 u d\mu + Ric_{m,n}(L)(\nabla f, \nabla f)u\right)d\mu$$
$$- \frac{2}{m-n}\int_M \tau\left(\nabla\phi.\nabla f + \frac{m-n}{2\tau}\right)^2 u d\mu$$

Corollary1.2 (X.-D.Li 2006) supposes that $Ric_{m,n}(L) \geq 0$ then $\tau \to \mu(\tau)$ is decreasing along the heat diffusion $(\partial_\tau - L)u = 0$.

## Perelman functional $\mathcal{F}$ and $\mathcal{W}$ for extend Ricci flow system

*We consider a system of evolution equations such that the stationary points satisfy*
$$Rc(g) = 2du \otimes du$$
$$\Delta^g u = 0$$
Bernhard list *in [3] extended the Ricci flow to the system*
$$\frac{\partial g(t)}{\partial t} = -2Rc(g(t)) + 4du \otimes du \qquad (1.14)$$
$$\frac{\partial u(t)}{\partial t} = \Delta^{g(t)} u(t)$$
*For a Riemannian metric $g(t)$, a function $u(t)$, and given initial data $g(0)$ and $u(0)$. this is a quasilinear, weakly parabolic, coupled system of second order. Here $du \otimes du$ is the tensor $\partial_i u \partial_j u dx^i \otimes dx^j$ and the laplacian of a function $u$ with respect to $g$ is given by $\Delta^g u = g^{ij}(\partial_i \partial_j u - \Gamma_{ij}^k \partial_k u)$.*

*Definition1.8(see[3]) let $\tau \in R$ be a positive real number. Then the entropy $\mathcal{W}$ of a configuration*
$$(g, u, f, \tau) \in \mathcal{M}(M) \times C^\infty(M) \times C^\infty(M) \times R^+$$
*is defined to be*
$$\mathcal{W}(g, u, f, \tau) := \int_M [\tau(S + |df|^2) + f - n](4\pi\tau)^{\frac{-n}{2}} e^{-f} dV$$
Where $S := R - 2|du|^2$ *(So the evolution of the metric can then be written as $\frac{\partial g_{ij}}{\partial t} = -2S_{ij}$) now we prove that $\mathcal{W}$ is scaling invariant.*

*Lemma1.7 (see[3]) let $\alpha > 0$ be a constant and $\varphi$ be a diffeomorphism of M. then the entropy $\mathcal{W}$ is invariant under simultaneous scaling of $g$ and $\tau$ by $\alpha$ in the sense that*
$$\mathcal{W}(\alpha g, u, f, \alpha \tau) = \mathcal{W}(g, u, f, \tau)$$
*and invariant under diffeomorphisms*
$$\mathcal{W}(\varphi^* g, \varphi^* u, \varphi^* f, \tau) = \mathcal{W}(g, u, f, \tau)$$
*Proof the invariance under diffeomorphisms is clear . also we have*

$$\mathcal{W}(\alpha g, u, f, \alpha \tau )$$
$$= \int_M [\alpha\tau(R(\alpha g) - 2(\alpha g)^{ij}\partial_i u\partial_j u + (\alpha g)^{ij}\partial_i f\partial_j f) + f$$
$$- n] (4\pi\alpha\tau)^{\frac{-n}{2}} e^{-f} \sqrt{\det(\alpha g)} dx$$
$$= \int_M [\alpha\tau(\alpha^{-1} R - 2\alpha^{-1}|du|^2 + \alpha^{-1}|df|^2) + f - n]\alpha^{\frac{-n}{2}} (4\pi\tau)^{\frac{-n}{2}} e^{-f} \alpha^{\frac{n}{2}} dV$$
$$= \mathcal{W}(g, u, f, \tau )$$

*Theorem1.15 (see[3])let M be a closed Riemannian manifold and assume that $g, u, f$ and $\tau$ satisfy on $[0, T) \times M$ the evolution equations*
$$\partial_t g = -2S_y$$
$$\partial_t u = \Delta u$$
$$\partial_t f = -\Delta f + |\nabla f|^2 - S + \frac{n}{2\tau}$$
$$\partial_t \tau = -1$$

*(S is the trace of symmetric tensor field $S_y$ )*
*then the following monotonicity formula holds:*
$$\partial_t \mathcal{W}(t) = \int_M \left(2\tau \left| S_y + \nabla^2 f - \frac{1}{2\tau} g\right|^2 + 4\tau|\nabla u - du(\nabla f)|^2\right) dm \geq 0$$

*Remark1.7(see[3]) the entropy $\mathcal{W}$ is non-decreasing and equality holds if and only if the solution is a homothetic shrinking gradient soliton. In this case $(g, u, f, \tau )(t)$ satisfies at every $t \in [0, T)$.*
$$S_y + \nabla^2 f - \frac{1}{2\tau} g = 0 \text{ and } \nabla u - du(\nabla f) = 0$$

*Theorem1.16(see[3]) if we vary $\mathcal{W}$ along the variation given by the following evolution equations*
$$\partial_t g = -2S_y - 2\nabla^2 f$$
$$\partial_t u = \Delta u - \langle du, df\rangle$$
$$\partial_t f = -\Delta f - S + \frac{n}{2\tau}$$
$$\partial_t \tau = -1$$

then we have

$$\partial_t \mathcal{W}(g, u, f, \tau )(t) = \int_M \left(2\tau \left| S_y + \nabla^2 f - \frac{1}{2\tau} g\right|^2 + 4\tau|\Delta u - \langle du, df\rangle|^2\right) dm \geq 0$$

*Definition 1.9(see[3])let $(g, u, \tau ) \in \mathcal{M}(M) \times C^\infty(M) \times R^+$ be a configuration .then we define*
$$\mu := \mu(g, u, \tau) = \inf_{f \in C^\infty(M)} \left\{ \mathcal{W}(g, u, f, \tau)| \int_M (4\pi\tau)^{\frac{-n}{2}} e^{-f} dV = 1 \right\}$$

*We investigate the remaining case of shrinking breathers.*

*Theorem1.17. (see[3])suppose $(g, u)(t)$ is a solution to (1.14) on $[0, T) \times M$ where M is closed. Fix a $\bar{\tau} \in [0, T)$ and define $\tau (t) := \bar{\tau} - t$. then $\mu(g, u, \tau )(t)$ is non-decreasing in t. If $\frac{\partial}{\partial t}\mu(t) = 0$ the solution is a gradient shrinking soliton.*

*Proposition1.4(see[3]) Let $(g, u)(t)$ be a shrinking breather on a closed manifold M. Then it necessarily is a gradient shrinking soliton.*

# Perelman's Functional $\mathcal{F}$ and $\mathcal{W}$ on Ricci Yang-Mills flow

*The Yang-Mills heat flow was first used by Atiyah (in [10]) and Bott and simon Donaldson. Donaldson used it to give an analytic proof of a Theorem of Narasimhan and Seshadri. also Atiyah and Bott used the Yang-Mills heat flow to study the topology of minimal Yang Mills connections. The Yang-Mills heat flow is a gauge-theoretic heat flow; that is, it is a differential equation for a field on a principal fiber bundle. In the study of geometric evolution equations monotonic quantities have always played an important role. Here we give a Review of the Ricci Yang- Mills flow using energy functional, but as for Ricci flow the resulting equations are not a-priori gradient equations. In other words we would like to follow the ideas of Perelman in order to write the Ricci Yang-Mills flow as a gradient flow. We claim that our coupled system is the gradient flow of some functional $\mathcal{F}(g, a, f)$ analogous to that of Perelman.*

*For definition of the Ricci Yang-Mills flow at first, we recall the $U(m)$ −vector bundle.*

*Let G be a unitary group $U(m) \subseteq GL(m, R) \subseteq GL(2m, R)$, a $G$ −vector bundle is a complex vector bundle of rank m together with Hermitian metric, a smooth function which assigns to each $p \in M$ a map*

$$\langle , \rangle_p : E_p \times E_p \to C$$

Which satisies the axioms
1. $\langle v, w \rangle_p$ is complex linear in $w$ and conjugate linear in $v$.
2. $\langle v, w \rangle_p = \overline{\langle w, v \rangle_p}$
3. $\langle v, w \rangle_p \geq 0$, with equality holding only if $v = 0$

*Definition 1.10 (see [9]) Let P be a $U(1)$ −bundle over a compact manifold. One can choose a metric on P such that the Ricci flow equations, with the additional hypothesis that size of the fiber remains fixed, yield the Ricci Yang-Mills flow:*

$$\frac{\partial g}{\partial t} = -2Rc + F^2$$

$$\frac{\partial A}{\partial t} = -D_A^* F(A)$$

*here "A" is a connection on P and $D^*$ is adjoint of the exterior derivative D. also F is a two-form on M where $F(A) = D_A A$ that DA denote the exterior covariant derivative. Recall that if G is a lie sub group of $GL(m, R)$, a $G$ −vector bundle is a rank m vector bundle whose transition functions take their values in G.*

*Definition 1.11 (see [9]) Let $(M, g)$ be a Riemannian manifold and let $E \to M$ denote a principal $K$ −bundle (K is a lie group) over M with connection A. In this section $\nabla$ will always refer to the Levi-Civita connection of g. Consider the functional*

$$\mathcal{F}(g, A, f) = \int_M \left( R - \frac{1}{4} |F|^2 + |\nabla f|^2 \right) e^{-f} dV$$

where R is the scalar curvature of the base metric, $f \in C^\infty(M)$ and $dV$ denotes the volume form of $g$.

We use the notation $\delta$ to refer to the first variation at 0 of other quantities with respect to the parameter $t$.

*Lemma1.8(see [9])* Let $\delta g_{ij} = v_{ij}$, $\delta A_i = \alpha_i$, $\delta f = h$, then

$$\delta \mathcal{F}(v, \alpha, h) = \int_M e^{-f} \left[ -v_{ij} \left( Rc_{ij} - \frac{1}{2}\eta_{ij} + \nabla_i \nabla_j f \right) - \alpha_j (\alpha^* F_j - \nabla^i f F_{ij}) \right.$$
$$\left. + \left(\frac{v}{2} - h\right)\left(2\Delta f - |\nabla f|^2 + R - \frac{1}{4}|F|^2\right) \right] dV$$

where $\mu_{ij} = F_i^k F_{kj}$

*Theorem1.18 (see [9])* Given $(g(t), A(t), f(t))$ a solution to following system

$$\frac{d}{dt}g_{ij} = -2Rc_{ij} + \eta_{ij} - 2\nabla_i \nabla_j f$$
$$\frac{d}{dt}A_i = -d^*F + \nabla^i f F_{ij} \qquad (1.15)$$
$$\frac{d}{dt}f = -\Delta f - R + \frac{1}{2}|F|^2$$

Then the functional $\mathcal{F}$ is monotonically increasing in $t$. In particular

$$\frac{d}{dt}\mathcal{F} = \int_M \left(2\left|Rc_{ij} - \frac{1}{2}\eta_{ij} + \nabla_i \nabla_j f\right|^2 + |d^*F - \nabla^i f F_{ij}|^2\right) e^{-f} dV \geq 0$$

*Definition1.12:* (we define a metric on our configuration space to be

$$\langle (g_1, A_1), (g_2, A_2) \rangle = \int (2(A_1, A_2) + 2(g_1, g_2)) e^{-f} dV$$

then the gradient flow of $\mathcal{F}$ becomes (1.15)

*Remark1.8 (see [9])* under equations (1.15) we know $\frac{d}{dt}\mathcal{F}(g(t), A(t), f(t)) \geq 0$. Equality is attained precisely when $R_{ij} - \frac{1}{2}F_i^k F_{kj} + \nabla_i \nabla_j f = 0$ and $d^*F_j - \nabla_i f F_{ij} = 0$ i.e. When $(g, A)$ is a steady gradient Ricci Yang-Mills soliton.

*Remark1.9* The solutions to equations (*1.15*) are equivalent to the system

$$\frac{d}{dt}g_{ij} = -2Rc_{ij} + F_{ik}F_{jk}$$
$$\frac{d}{dt}A_i = -d^*F_i \qquad (1.16)$$
$$\frac{d}{dt}f = -\Delta f + |\nabla f|^2 - R + \frac{1}{2}|F|^2$$

*corollary 1.3* under equations (1.16)

$$\frac{d}{dt}\mathcal{F}(g(t), A(t), f(t)) \geq 0$$

*Proposition1.5 (see[9])* there exists a unique minimizer $\bar{f}$ of $\mathcal{F}(g, a, f)$ subject to the constraint

$$\int e^{-f} dV = 1$$

*Definition 1.13* according to this proposition we can then define

$$\lambda(g, A) = \inf\left\{ \mathcal{F}(g, A, f) : f \in C^\infty, \int e^{-f} dV = 1 \right\}$$

*Proposition1.6 (see [9]) if $(g(.), A(.))$ is a solution to the Ricci Yang - Mills flow, then $\lambda(g(t), A(t))$ is non - decreasing in time.*

*Remark1.10 (see [9]) the minimum value of $\lambda(g, A)$ is equal to $\lambda_1(g, A)$, where $\lambda_1(g, A)$ is the smallest eigenvalue of the elliptic operator $-4\Delta + R - \frac{1}{4}|F|^2$.*

*Then the minimizer, $f_0$, of $\mathcal{F}$ satisfies the Euler - Lagrange equation*
$$\lambda(g, A) = 2\Delta f_0 - |\nabla f_0|^2 + R - \frac{1}{4}|F|^2$$

*Definition1.14 A solution $(g(t), A(t))$ to the Ricci Yang -Mills flow is called a breather if there exist times $t_1 < t_2$, a constant $\alpha$, and a diffeomorphism $\varphi \colon M \to M$ such that*
$$g(t_2) = \alpha \varphi^* g(t_1), A(t_2) = \alpha \varphi^* A(t_1)$$
*$\alpha > 1, \alpha < 1$ and $\alpha = 1$ correspond to $(g(t), A(t))$ being a expanding, shrinking or steady breather respectively.*

*Theorem1.19(see[9]) let $(M^n, g(t), A(t))$ be a solution to the Ricci Yang-Mills flow on a closed manifold. if there exist $t_1 < t_2$ with $\lambda(g(t_1), A(t_1)) = \lambda(g(t_2), A(t_2))$ then $(g(t), A(t))$ is a steady gradient Ricci Yang-Mills soliton, which must have $|F|^2 = 0$ and be scalar flat.*

*Definition 1.15(see[9])we define*
$$\mathcal{W}(g, A, f, \tau) = \int_M \left( \tau \left( |\nabla f|^2 + R + \frac{1}{4}|F|^2 \right) + f - n \right) (4\pi\tau)^{\frac{-n}{2}} e^{-f} dV$$
*where $\tau = T - t$ and $(g(t), A(t))$ is a solution to Ricci Yang - Mills flow which exists on a maximal time interval of the form $[0, T]$ where $T < \infty$.*

*Theorem1.20(see[9]) let $v_{ij} = \delta g_{ij}$, $\delta A_i = \alpha_i$, $\delta f = h$ and $\delta \tau = \sigma$. Then*
$$\delta \mathcal{W}(v, \alpha, h, \sigma) = \int_M (4\pi\tau)^{\frac{-n}{2}} e^{-f} dV \left[ \sigma \left( |\nabla f|^2 + R + \frac{1}{4}|F|^2 \right) - \tau v_{ij} \left( Rc_{ij} + \frac{1}{2}\eta_{ij} + \nabla_i \nabla_j f \right) \right.$$
$$- \tau \alpha_j \left( d^* F_j - \nabla^i f F_{ij} \right) + h$$
$$\left. + \left[ \tau \left( 2\Delta f - |\nabla f|^2 + R - \frac{1}{4}|F|^2 \right) + f - n \right] \left( \frac{v}{2} - h - \frac{n\sigma}{2\tau} \right) \right]$$

*Remark 1.11(see[9]) Consider the following system of equations*
$$\frac{d}{dt} g_{ij} = -2 \left( Rc_{ij} - \frac{1}{2}\eta_{ij} + \nabla_i \nabla_j f \right)$$
$$\frac{d}{dt} A_i = -d^* F_i + \nabla^i f F_{ij}$$
$$\frac{d}{dt} f = -\Delta f - R + \frac{1}{2}|F|^2 + \frac{n}{2\tau}$$
$$\frac{d}{dt} \tau = -1$$
*for a solution to this system we have*
$$\frac{d\mathcal{W}}{dt} = \int_M \left[ 2\tau \left| Rc_{ij} - \frac{1}{2\tau} g_{ij} + \nabla_i \nabla_j f \right|^2 + \tau \left| d^* F_j + \nabla^i f F_{ij} \right|^2 + \frac{1}{4}|F|^2 \right.$$
$$\left. - \frac{1}{2} \tau |\eta|^2 \right] (4\pi\tau)^{\frac{-n}{2}} e^{-f} dV$$

*Definition1.15(see[9]) let $(M, g(t), A(t))$ be a solution to RYM-flow which exists on a maximal time interval $T < \infty$. $(M, g(t), A(t))$ is a low - energy solution if*
$$\lim_{t \to T}(T - t)\, |F|^2_{C^0(M_t)} = 0$$
*So according to this Definition we get the following corollary.*

*Corollary1.4 (see[9]) Let $(M, g(t), A(t))$ be a low energy solution to RYM flow on $[0, T)$ then there exists $t_0 < T$ such that for all $t_0 \leq t < T$, we have*
$$\frac{d\mathcal{W}}{dt} \geq 0$$

*Definition 1.16 Given $(M, g, A)$ at $t \in R$ let*
$$\mu(g, A, \tau) = \inf_f \left\{ \mathcal{W}(g, A, f, \tau) \Big| \int_M (4\pi\tau)^{\frac{-n}{2}} e^{-f} dV = 1 \right\}$$

*Corollary (see[9]) let $(M, g(t), A(t))$ be a low - energy solution to RYM flow on $[0, T)$. Then there exists $t_0 < T$ such that for all $t_0 \leq t < T$ we have $\frac{d\mu}{dt} \geq 0$.*

## 2. Perelman Reduced volume and reduced length

*In this section we discuss Perelman's notions of the $\mathcal{L}$ —length in the context of Ricci flows. This is a functional defined on paths in space –time parametrized by backward time, denoted $\tau$. The main purpose of this section is to use the Li-Yau-Perelman distance to define the Perelman's reduced volume, which was introduced by Perelman, and prove the monotonicity property of the reduced volume under the Ricci flow. The reduced distance, i.e. the $l$ —function. The $l$ —function is defined in terms of a natural curve energy along the Ricci flow, which is analogous to the classical curve energy employed in the study of geodesics, but involves the evolving metric, as well as the scalar curvature as a potential term. There are two applications of this theory of Perelman. We use the theory of $\mathcal{L}$ —geodesics and the associated notion of reduced volume to establish non-collapsing results. The second application will be to $\kappa$ —non-collapsed solutions of bounded non-negative curvature.*

## Perelman reduced distance

*Definition2.1: suppose that either M is compact or $g_{ij}(\tau)$ are complete and have uniformly bounded curvature. The $\mathcal{L}$ —length of a smooth space curve $\gamma: [\tau_1, \tau_2] \to M$ is defined by*

$$\mathcal{L}(\gamma) = \int_{\tau_1}^{\tau_2} \sqrt{\tau}\left(R(\gamma(\tau)) + |\dot{\gamma}(\tau)|^2\right) d\tau \quad (2.1)$$

*Where $\tau = \tau(t)$ satisfying $\frac{d\tau}{dt} = -1$ and the scalar curvature $R(\gamma(\tau))$ and the norm $|\dot{\gamma}(\tau)|$ are evaluated using the metric at time $t = t_0 - \tau$.*

We consider an 1-parameter family of curves $\gamma_s : [\tau_1, \tau_2] \to M$ parameterized by $s \in (-\varepsilon, \varepsilon)$. Equivalently, We have a map $\tilde{\gamma}(s, \tau)$ with $s \in (-\varepsilon, \varepsilon)$ and $\tau \in [\tau_1, \tau_2]$. Putting $X = \frac{\partial \tilde{\gamma}}{\partial \tau}$ and $Y = \frac{\partial \tilde{\gamma}}{\partial s}$, We have $[X, Y] = 0$ (Because $\nabla_X Y - \nabla_Y X - [X, Y] = 0$). This implies that $\nabla_X Y = \nabla_Y X$.

Writing $\delta_Y$ as shorthand for $\frac{d}{ds}\big|_{s=0}$, and restricting to the curve $\gamma(\tau) = \tilde{\gamma}(0, \tau)$, we have $(\delta_Y Y)(\tau) = Y(\tau)$ and $(\delta_Y X)(\tau) = (\nabla_X Y)(\tau)$. Set $s = \sqrt{\tau}$ one sees immediately with respect to the variable $s$.

The $\mathcal{L}-$functional is

$$\mathcal{L}(\gamma) = \int_{\sqrt{\tau_1}}^{\sqrt{\tau_2}} \left(\frac{1}{2}\left|\frac{d\gamma}{ds}\right|^2 + 2s^2 R(\gamma(s))\right) ds \quad (2.2)$$

Remark2.2($\mathcal{L}$ on Riemannian Products) suppose that we are given a Riemannian product solution $\left(N_1^{n_1} \times N_2^{n_2}, h_1(\tau) + h_2(\tau)\right)$ to the backward Ricci flow and a $C^1-$Path $\gamma = (\alpha, \beta) : [\tau_1, \tau_2] \to N_1 \times N_2$ so $\mathcal{L}_{h_1+h_2}(\gamma) = \mathcal{L}_{h_1}(\alpha) + \mathcal{L}_{h_2}(\beta)$.

Example suppose that Ricci flow is a constant family of Euclidean metrics on $R^n \times [0, T]$. Then we have $R(\gamma(\tau)) = 0$. So according to (2.2) we get

$$\mathcal{L}(\gamma) = \frac{1}{2}\int_0^{\bar{\tau}} \sqrt{\tau}\left|\frac{d\gamma}{d\tau}\right|^2 d\tau$$

If we change $= \sqrt{\tau}$, $\frac{d\gamma}{d\tau} = \frac{1}{2\sqrt{\tau}}\frac{d\gamma}{ds}$, $d\tau = 2sds$

$$\mathcal{L}(\gamma) = \frac{1}{2}\int_0^{\bar{s}^2} s\frac{1}{4s^2}\left|\frac{d\gamma}{ds}\right|^2 2sds = \frac{1}{2}\int_0^{\bar{s}^2}\left|\frac{d\gamma}{ds}\right|^2 ds$$

So we arrive to (2.2) formula.

**First variation formulae for $\mathcal{L}-$geodesics**

At first recall the classical variation formulae that give the first derivative of the metric function $d(x_0, x)$ on a Riemannian manifold $(M, g)$.

Let $\gamma : [0,1] \to M$ ranges over all $C^1-$Curves from $x_0$ to $x$, and the Dirichlet energy $E(\gamma)$ of the curve is given by the formula

$$E(\gamma) = \frac{1}{2}\int_0^1 |X|_g^2 dt \quad \text{where } X = \frac{\partial \gamma}{\partial t}$$

We know that the distance $d(x_0, x)$ on a Riemannian manifold $(M, g)$ can be defined by the energy-minimisation formula

$$\tfrac{1}{2}d(x_0, x)^2 = \inf_\gamma E(\gamma) \quad (2.3)$$

Also if $\gamma: [0, \bar{t}] \to M$, $\gamma(0) = p$, $\gamma(\bar{t}) = q$

$$inf\{E(\gamma) | \gamma: [0, \bar{t}] \to M, \gamma(0) = p, \gamma(\bar{t}) = q\} = \frac{d^2(p,q)}{\bar{t}}$$

Lemma 2.1 The first variation formula for Dirichlet energy $E(\gamma)$ is

$$\frac{d}{ds}E(\gamma) = \langle X, Y \rangle|_{t=0}^{1} - \int_0^1 \langle X, \nabla_Y Y \rangle dt \quad (2.4)$$

Proof: consider $D_t X$ means $D_t X = \nabla_{\frac{\partial}{\partial t}} X$ and $X = \frac{\partial \gamma}{\partial s}$, $\dot{\gamma} = \frac{\partial \gamma}{\partial t}$.

$$\delta E(\gamma) = \delta \frac{1}{2} \int_0^1 \langle \dot{\gamma}, \dot{\gamma} \rangle dt$$

$$= \frac{1}{2} \frac{\partial}{\partial s}\bigg|_{s=0} \int_0^1 \langle \frac{\partial \gamma}{\partial t}, \frac{\partial \gamma}{\partial t} \rangle dt$$

$$= \frac{1}{2} \int_0^1 \langle D_s \left(\frac{\partial \gamma}{\partial t}\right), \frac{\partial \gamma}{\partial t} \rangle + \langle \frac{\partial \gamma}{\partial t}, D_s \left(\frac{\partial \gamma}{\partial t}\right) \rangle dt$$

$$= \int_0^1 \langle D_s \left(\frac{\partial \gamma}{\partial t}\right), \frac{\partial \gamma}{\partial t} \rangle dt$$

$$= \int_0^1 \langle D_t X, \dot{\gamma} \rangle dt$$

With integration by parts we get

$$= \langle X, \dot{\gamma} \rangle|_0^1 - \int_0^1 \langle X, D_t \dot{\gamma} \rangle dt$$

$$= \langle X, Y \rangle|_0^1 - \int_0^1 \langle X, \nabla_Y Y \rangle$$

Remark 2.3 if we fix the end points of $\gamma$, then the first term of the right hand side of (2.4) vanishes. If we consider arbitrary variations $X$ of $\gamma$ with fixed end points, we thus conclude that in order to be a minimiser for (2.3), that $\gamma$ must obey the geodesic flow equation

$$\nabla_Y Y = 0$$

Now we develop analogous variational formulae for $\mathcal{L}$ −length (reduced distance) on Ricci flow

Theorem 2.1 (Perelman [5]) the first variation formula for $\mathcal{L}$ −length is

$$\delta_Y(\mathcal{L}) = 2\sqrt{\tau}\langle X, Y \rangle|_{\tau_1}^{\tau_2} + \int_{\tau_1}^{\tau_2} \sqrt{\tau} \langle Y, \nabla R - 2\nabla_X X - 4Ric(.,X) - \frac{1}{\tau}X \rangle d\tau \quad (2.5)$$

Where $\langle .,. \rangle$ denotes the inner product with respect to the metric $g_{ij}(\tau)$.

Proof let $X(\tau) = \dot{\gamma}(\tau)$, so according to (2.1) we obtain

$$\delta_Y(\mathcal{L}) = \int_{\tau_1}^{\tau_2} \sqrt{\tau}(\nabla_Y R + 2\langle \nabla_X Y, X\rangle)d\tau$$

$$= \int_{\tau_1}^{\tau_2} \sqrt{\tau}(\langle \nabla R, Y\rangle + 2\langle \nabla_X Y, X\rangle)d\tau$$

Because $\frac{\partial \tau}{\partial t} = -1$ so $\frac{dg_{ij}}{d\tau} = 2R_{ij}$ ,using this fact we get

$$\frac{d\langle Y,X\rangle}{d\tau} = \langle \nabla_X Y, X\rangle + \langle Y, \nabla_X X\rangle + 2Ric(Y,X) \quad (2.6)$$

Because we can break $\frac{d\langle Y,X\rangle}{d\tau}$ into two parts :first assumes that the metric is constant and the second deals with the variation with $\tau$ of the metric .The first contribution is the usual formula

$$\frac{d\langle Y,X\rangle}{d\tau} = \langle \nabla_X Y, X\rangle + \langle Y, \nabla_X X\rangle$$

We show last term in equation (2.6) is that come from differentiating the metric with respect to $\tau$ and Ricci flow equation we get

$$\frac{d\langle Y,X\rangle}{d\tau} = 2Ric(Y,X)$$

So we arrive to equality of (2.6)

We continue the proof of (2.5)

$$\int_{\tau_1}^{\tau_2} \sqrt{\tau}(\langle Y, \nabla R\rangle + 2\langle \nabla_X Y, X\rangle)d\tau =$$

$$= \int_{\tau_1}^{\tau_2} \sqrt{\tau}\left(\langle Y, \nabla R\rangle + 2\frac{d\langle Y,X\rangle}{d\tau} - 2\langle Y, \nabla_X X\rangle - 4Ric(Y,X)\right)d\tau$$

$$= 2\sqrt{\tau}\langle X, Y\rangle\Big|_{\tau_1}^{\tau_2} + \int_{\tau_1}^{\tau_2} \sqrt{\tau}\,\langle Y, \nabla R - 2\nabla_X X - 4Ric(X,.) - \frac{1}{\tau}X\rangle d\tau$$

Lemma2.2 (Perelman) we consider a variation $\gamma(\tau, s)$ with fixed endpoints, so that $Y(\tau_1) = Y(\tau_2) = 0$. so from (2.5) the $\mathcal{L}$ −shortest geodesic $\gamma(\tau, s)$ for $\tau \in [\tau_1, \tau_2]$ satisfies the following $\mathcal{L}$ − geodesic equation a

$$\nabla_X X - \frac{1}{2}\nabla R + \frac{1}{2\tau}X + 2Ric(X,.) = 0 \quad (2.7)$$

This equation is called Euler-Lagrange Equation.

Therefore we can say for any $\tau_2 > \tau_1 > 0$ ,there is always an $\mathcal{L}$ −geodesic $\gamma(\tau)$ for $\tau \in [\tau_1, \tau_2]$ .

Written with respect to the variable $s = \sqrt{\tau}$ the $\mathcal{L}$ −geodesic equation becomes

$$\nabla_{\hat{X}}\hat{X} - 2s^2\nabla R + 4sRic(\hat{X},.) = 0$$

Where $\hat{X} = \frac{d\gamma}{ds} = 2sX$

**Lemma2.3(Perelman)** let $\gamma: [\tau_1, \tau_2] \to M$ be an $\mathcal{L}$-geodesic. Then $\lim_{\tau \to 0} \sqrt{\tau} X(\tau)$ exists.

*Proof:* multiplying the $\mathcal{L}$-geodesic equation (2.7) by $\sqrt{\tau}$, we get

$$\nabla_X(\sqrt{\tau}\, X) = \frac{\sqrt{\tau}}{2} \nabla R - 2\sqrt{\tau}\, Ric(X,.) \quad On \quad [\tau_1, \tau_2]$$

*or equivalently*

$$\frac{d}{d\tau}(\sqrt{\tau}\, X) = \frac{\sqrt{\tau}}{2} \nabla R - 2\, Ric(\sqrt{\tau} X,.) \quad On \quad [\tau_1, \tau_2]$$

*Thus if a continuous curve, defined on $[\tau_1, \tau_2]$, satisfies the $\mathcal{L}$-geodesic equation on every subinterval $0 < \tau_1 \leq \tau \leq \tau_2$, then $\sqrt{\tau_1} X(\tau_1)$ has a limit as $\tau_1 \to 0^+$.*

*Remark2.4 for a fixed $p \in M$, by taking $\tau_1 = 0$ and $\gamma(0) = p$, the vector $v = \lim_{\tau \to 0} \sqrt{\tau} X(\tau)$ is well-defined in $T_P M$. The $\mathcal{L}-$ exponential map $\mathcal{L}exp_\tau: T_P M \to M$ sends $v$ to $\gamma(\tau)$.*

*Note that for any vector $v \in T_P M$, we can find an $\mathcal{L}$-geodesic $\gamma(\tau)$ with $\lim_{\tau \to 0^+} \sqrt{\tau} \dot{\gamma}(\tau) = v$.*

*Now we give an Estimate for speed of $\mathcal{L}$-geodesics*

*Theorem(see[4])2.2:* let $(M^n, g(\tau))$, $\tau \in [0, T]$, be a solution to the backward Ricci flow with bounded sectional curvature and $\max\{|Rm|, |Rc|\} \leq C_0 < \infty$ on $M \times [0, T]$. There exists a constant $C(n) < \infty$ depending only on $n$ such that given $0 \leq \tau_1 \leq \tau_2 < T$, if $\gamma: [\tau_1, \tau_2] \to M$ is an $\mathcal{L}$-geodesic with

$$\lim_{\tau \to \tau_1} \sqrt{\tau} \frac{d\gamma}{d\tau}(\tau) = v \in T_{\gamma(\tau_1)} M$$

*Then for any $\in [\tau_1, \tau_2]$,*

$$\tau \left| \frac{d\gamma}{d\tau}(\tau) \right|^2_{g(\tau)} \leq e^{6C_0 T} |v|^2 + \frac{C(n) T}{\min\{T - \tau_2, C_0^{-1}\}} (e^{6C_0 T} - 1)$$

Where $|v|^2 = |v|^2_{g(\tau_1)}$.

*Also we give a lemma about existence of $\mathcal{L}$-geodesics between any two space-time endpoints.*

*Lemma2.4(see[5])* let $(M^n, g(\tau))$, $\tau \in [0, T]$, be a complete solution to the backward Ricci flow with bounded sectional curvature. Given $p, q \in M$ and $0 < \tau_1 \leq \tau_2 < T$, There exists a smooth path $\gamma(\tau): [\tau_1, \tau_2] \to M$ from $p$ to $q$ such that $\gamma$ has the minimal of $\mathcal{L}$-length among all such paths. Furthermore, all $\mathcal{L}$-length minimizing paths are smooth $\mathcal{L}$-geodesics.

*Corollary2.1 (see[8]):* if we extend the curve $\gamma$ for piecewise smooth curves (where $0 \leq \tau_1 < \tau_2 \leq t_0$) then for first variation formula of breaking points $\tau_0' = \tau_1 < \tau_2' < \cdots < \tau_k' = \tau_2$, we have

$$\delta(\mathcal{L}(\gamma)) =$$
$$\int_{\tau_1}^{\tau_2} \sqrt{\tau} \langle \nabla R - 2\nabla_X X - 4Ric(X) - \frac{1}{\tau}X, Y \rangle d\tau + 2\sqrt{\tau}\langle X, Y\rangle_\tau \Big|_{\tau_1}^{\tau_2} + \sum_{i=2}^{k-1} 2\sqrt{\tau_i}\langle X^-(\tau_i{}') - X^+(\tau_i{}'), Y(\tau_i{}')\rangle_{\tau_i{}'}$$

*Definition 2.2: Fixing a point $p$, we denote by $L(q, \bar{\tau})$ the $\mathcal{L}$-length of the $\mathcal{L}$-shortest curve $\gamma(\tau), 0 \leq \tau \leq \bar{\tau}$, joining $p$ and $q$. In other words*

$$L(q,\bar{\tau}) = \inf\{\mathcal{L}(\gamma) | \gamma: [0,\bar{\tau}] \to M \text{ with } \gamma(0) = p, \gamma(\bar{\tau}) = q\}$$

*Theorem 2.3 (Perelman[5]) suppose that*

$$H(X) = -\frac{\partial R}{\partial \tau} - \frac{R}{\tau} - 2\langle \nabla R, X\rangle + 2Ric(X,X)$$

And

$$K = \int_0^{\bar{\tau}} \tau^{\frac{3}{2}} H(X) d\tau$$

Then

a) $|\nabla L|^2 = -4\bar{\tau}R + \frac{2}{\sqrt{\bar{\tau}}}L - \frac{4}{\sqrt{\bar{\tau}}}K$
b) $\frac{\partial L}{\partial \bar{\tau}} = 2\sqrt{\bar{\tau}}R - \frac{1}{2\bar{\tau}}L + \frac{1}{\bar{\tau}}K$

*Proof: The first variation formula in (2.5) implies that*

$$\nabla_Y L(q, \bar{\tau}) = \langle 2\sqrt{\bar{\tau}}X(\bar{\tau}), Y(\bar{\tau})\rangle$$

So $\nabla L(q,\bar{\tau}) = 2\sqrt{\bar{\tau}}X(\bar{\tau})$

*At first we prove*

$$\frac{\partial L}{\partial \bar{\tau}}(\gamma(\bar{\tau}), \bar{\tau}) = \frac{d}{d\tau}L(\gamma(\tau), \tau)\Big|_{\tau=\bar{\tau}} - \langle \nabla L, X\rangle \quad (2.8)$$

Since $\gamma(\bar{\tau}) = -\frac{\partial}{\partial t} + X(\bar{\tau})$, So the chain rule implies

$$\frac{d}{d\tau}L(\gamma(\tau), \tau)\Big|_{\tau=\bar{\tau}} = \frac{d}{d\bar{\tau}}L(\gamma(\bar{\tau}), \bar{\tau}) + \langle \nabla L, X\rangle$$

*Therefore we get to (2.8)*

*Also from (2.1) and (2.8) we obtain*

$$\frac{\partial L}{\partial \bar{\tau}}(\gamma(\bar{\tau}), \bar{\tau}) = \sqrt{\bar{\tau}}(R + |X|^2) - \langle \nabla L, X\rangle$$

*And because $\nabla L(q,\bar{\tau}) = 2\sqrt{\bar{\tau}}X(\bar{\tau})$ so*

$$\frac{\partial L}{\partial \bar{\tau}}(\gamma(\bar{\tau}), \bar{\tau}) = \sqrt{\bar{\tau}}(R + |X|^2) - \sqrt{\bar{\tau}}|X|^2$$

$$= 2\sqrt{\bar{\tau}}R - \sqrt{\bar{\tau}}(R + |X|^2) \quad (2.9)$$

Now we compute the $+|X|^2$.

According to (2.5) and this fact that $\gamma'(\tau) = \frac{\partial}{\partial \tau} + X(\tau)$, we get

$$\frac{d}{d\tau}(R(\gamma(\tau), \tau) + |X(\tau)|^2) = \frac{\partial R}{\partial \tau} + \langle \nabla R, X \rangle + |X(\tau)|^2$$

Also we know $|X(\tau)|^2 = \langle X(\tau), X(\tau) \rangle$

So from $\frac{d}{d\tau}|X(\tau)|^2 = \frac{d}{d\tau}\langle X(\tau), X(\tau) \rangle = 2\langle \nabla_X X, X \rangle + 2Ric(X, X)$

We get $\frac{d}{d\tau}(R(\gamma(\tau), \tau) + |X(\tau)|^2) = \frac{\partial R}{\partial t} + \langle \nabla R, X \rangle + 2\langle \nabla_X X, X \rangle + 2Ric(X, X)$

Also we have $\nabla_X X = \frac{1}{2}\nabla R - \frac{1}{2\tau}X - 2Ric(X,.)$ so

$$\frac{d}{d\tau}(R(\gamma(\tau), \tau) + |X(\tau)|^2) = \frac{\partial R}{\partial t} + \frac{1}{\tau}R + 2\langle \nabla R, X \rangle - 2Ric(X, X) - \frac{1}{\tau}(R + |X|^2)$$

$$= -H(X) - \frac{1}{\tau}(R + |X|^2)$$

So $\frac{d}{d\tau}\left(\tau^{\frac{3}{2}}(R + |X|^2)\right)\bigg|_{\tau=\bar{\tau}} = \frac{1}{2}\sqrt{\bar{\tau}}(R + |X|^2) - \bar{\tau}^{\frac{3}{2}}H(X)$

$$= \frac{1}{2}\frac{d}{d\tau}L(\gamma(\tau), \tau)\bigg|_{\tau=\bar{\tau}} - \bar{\tau}^{\frac{3}{2}}H(X)$$

So $\bar{\tau}^{\frac{3}{2}}(R + |X|^2) = \frac{1}{2}L(q, \bar{\tau}) - K$ where

$$K = \int_0^{\bar{\tau}} \tau^{\frac{3}{2}}H(X)d\tau$$

Therefore by applying $|\nabla L|^2 = 4\bar{\tau}|X|^2 = -4\bar{\tau}R + 4\bar{\tau}(R + |X|^2)$ we conclude a and b. So proof is complete ∎

**Second variation formulae for $\mathcal{L}$ −geodesics**

*Proposition 2.1 If we compute the second variation of energy when $\gamma$ is a geodesic, we get*

$$\frac{d^2}{ds^2}E(\gamma) = \langle \nabla_Y Y, \gamma \rangle|_0^\tau + \int_0^\tau (\langle \nabla_Y \nabla_Y X, X \rangle + \langle \nabla_Y X, \nabla_Y X \rangle)dt$$

*Now turn to the second spatial variation of the reduced length.*

*Theorem 2.4(Perelman [5]) for any $\mathcal{L}-$geodesic $\gamma$, we have*

$\delta_Y^2(\mathcal{L}) = 2\sqrt{\tau}\langle\nabla_Y Y, X\rangle|_0^{\bar{\tau}} + \int_0^{\bar{\tau}} \sqrt{\tau}[2|\nabla_X Y|^2 + 2\langle R(Y,X)Y, X\rangle + \nabla_Y\nabla_Y R + 2\nabla_X Ric(Y,Y) - 4\nabla_Y Ric(Y,X)]d\tau$.

*Proof to compute the second variation $\delta_Y^2(\mathcal{L})$. we start the first variation formula.*

*We know*

$$\delta_Y(\mathcal{L}) = \int_0^{\bar{\tau}} \sqrt{\tau}(\nabla_Y R + 2\langle\nabla_X Y, X\rangle)d\tau$$

so $\qquad \delta_Y^2(\mathcal{L}) = Y\left(\int_0^{\bar{\tau}} \sqrt{\tau}(\nabla_Y R + 2\langle\nabla_X Y, X\rangle)d\tau\right)$

*We know* $\nabla_Y R = YR$, *so* $\delta_Y^2(\mathcal{L}) = Y\left(\int_0^{\bar{\tau}} \sqrt{\tau}(YR + 2\langle\nabla_X Y, X\rangle)d\tau\right)$

$$= \left(\int_0^{\bar{\tau}} \sqrt{\tau}\big(Y(Y(R)) + 2\langle\nabla_Y\nabla_X Y, X\rangle + 2|\nabla_X Y|^2\big)d\tau\right)$$

*Because $\nabla_X Y = \nabla_Y X$ we get*

$$= \left(\int_0^{\bar{\tau}} \sqrt{\tau}\big(Y(Y(R)) + 2\langle\nabla_Y\nabla_X Y, X\rangle + 2|\nabla_X Y|^2\big)d\tau\right)$$

*Also we know $\nabla_X\nabla_Y Z - \nabla_Y\nabla_X Z - \nabla_{[X,Y]}Z = R(X,Y)Z$ for all vectors $X,Y,Z$. so*

$$2\langle\nabla_Y\nabla_X Y, X\rangle = 2\langle\nabla_X\nabla_Y Y, X\rangle + 2\langle R(Y,X)Y, X\rangle$$

*For continuing proof, we prove the following equality*

$\frac{\partial}{\partial\tau}\langle\nabla_Y Y, X\rangle = \langle\nabla_X\nabla_Y Y, X\rangle + \langle\nabla_Y Y, \nabla_X X\rangle + 2Ric(\nabla_Y Y, X) + \langle\frac{\partial}{\partial\tau}\nabla_Y Y, X\rangle$ $\quad$ (2.10)

*We can break $\frac{\partial}{\partial\tau}\langle\nabla_Y Y, X\rangle$ into two parts: the first assumes that the metric is constant and the second deals with the variation with $\tau$ of the metric. The first contribution is the usual formula*

$$\frac{\partial}{\partial\tau}\langle\nabla_Y Y, X\rangle = \langle\nabla_X\nabla_Y Y, X\rangle + \langle\nabla_Y Y, \nabla_X X\rangle$$

*This gives the first two terms of the right-hand side of the equation. we show that the last two terms in that equation come from differentiating the metric with respect to $\tau$. To do this recall that in local coordinates, writing the metric as $g_{ij}$, we have*

$$\langle\nabla_Y Y, X\rangle = g_{ij}\big(Y^k\partial_k Y^i + \Gamma_{kl}^i Y^k Y^l\big)X^j$$

*There are two contributions coming from differentiating the metric with respect to $\tau$. The first is when we differentiate $g_{ij}$. This leads to*

$$2Ric_{ij}\big(Y^k\partial_k Y^i + \Gamma_{kl}^i Y^k Y^l\big)X^j = 2Ric(\nabla_Y Y, X)$$

*The other contribution is from differentiating the Christoffel symbols .This yields*

$$g_{ij}\frac{\partial \Gamma^i_{kl}}{\partial \tau}Y^kY^lX^j$$

*Differentiating the formula* $\Gamma^i_{kl} = \frac{1}{2}g^{si}(\partial_k g_{sl} + \partial_l g_{sk} - \partial_s g_{kl})$ *leads to*

$$g_{ij}\frac{\partial \Gamma^i_{kl}}{\partial \tau} = -2Ric_{ij}\Gamma^i_{kl} + g_{ij}g^{si}(\partial_k Ric_{sl} + \partial_l Ric_{sk} - \partial_s Ric_{kl})$$

$$= -2Ric_{ij}\Gamma^i_{kl} + \partial_k Ric_{jl} + \partial_l Ric_{jk} - \partial_j Ric_{kl}$$

*Thus, we have*

$$g_{ij}\frac{\partial \Gamma^i_{kl}}{\partial \tau}Y^kY^lX^j = \left(-2Ric_{ij}\Gamma^i_{kl} + \partial_k Ric_{jl}\right)Y^kY^lX^j$$

$$= \nabla_X Ric(Y,Y)$$

*so we obtain*

$$2\langle \nabla_Y \nabla_X Y, X\rangle = 2\frac{d}{d\tau}\langle \nabla_Y Y, X\rangle - 4Ric(\nabla_Y Y, X) - 2\langle \nabla_Y Y, \nabla_X X\rangle - 2\langle \frac{\partial}{\partial \tau}\nabla_Y Y, X\rangle + 2\langle \nabla_X \nabla_Y Y, X\rangle + 2\langle R(Y,X)Y, X\rangle$$

*also we can compute*

$$\langle \frac{\partial}{\partial \tau}\nabla_Y Y, X\rangle = 2(\nabla_Y Ric)(Y,X) - (\nabla_X Ric)(Y,Y) \quad (2.11)$$

*hence from (2.10) and (2.11)we get*

$$\frac{\partial}{\partial \tau}\langle \nabla_Y Y, X\rangle = \langle \nabla_X \nabla_Y Y, X\rangle + \langle \nabla_Y Y, \nabla_X X\rangle + 2Ric(\nabla_Y Y, X) + 2(\nabla_Y Ric)(Y,X) - (\nabla_X Ric)(Y,Y)$$
(2.12)

*Suppose* $Y(0) = 0$ *and the fact that* $\sqrt{\tau}X(\tau)$ *has a limit as* $\tau \to 0$ *are used to get the third equality below .So applying (2.12) and integrating by parts, we arrive*

$$\delta^2_Y(\mathcal{L}(\gamma)) = \left(\int_0^{\bar{\tau}}\sqrt{\tau}(Y(Y(R)) + 2\langle R(Y,X)Y,X\rangle + 2|\nabla_Y X|^2)d\tau\right) + 2\int_0^{\bar{\tau}}\sqrt{\tau}\Big(\frac{\partial}{\partial \tau}\langle \nabla_Y Y, X\rangle -$$

$$\langle \nabla_Y Y, \nabla_X X\rangle - 2Ric(\nabla_Y Y, X) - 2(\nabla_Y Ric)(Y,X) + (\nabla_X Ric)(Y,Y)\Big)d\tau$$

$$=$$
$$\left(\int_0^{\bar{\tau}}\sqrt{\tau}(Y(Y(R)) + 2\langle R(Y,X)Y,X\rangle + 2|\nabla_Y X|^2)d\tau\right) + 2\int_0^{\bar{\tau}}\sqrt{\tau}(-\langle \nabla_Y Y, \nabla_X X\rangle - 2Ric(\nabla_Y Y, X) -$$

$$2(\nabla_Y Ric)(Y,X) + (\nabla_X Ric)(Y,Y))d\tau + 2\sqrt{\tau}\langle \nabla_Y Y, X\rangle\Big|_0^{\bar{\tau}} - \int_0^{\bar{\tau}}\frac{1}{\sqrt{\tau}}\langle \nabla_Y Y, X\rangle d\tau = 2\sqrt{\bar{\tau}}\langle \nabla_Y Y, X\rangle +$$

$$\left(\int_0^{\bar{\tau}}\sqrt{\tau}(Y(Y(R)) - \nabla_Y Y.\nabla R + 2\langle R(Y,X)Y,X\rangle + 2|\nabla_Y X|^2)\right)d\tau + 2\int_0^{\bar{\tau}}\sqrt{\tau}\Big(-\langle \nabla_Y Y, \big[\nabla_X X +$$

$$2Ric(X) - \frac{1}{2}\nabla R + \frac{1}{2\tau}X]\rangle - 2\left((\nabla_Y Ric)(Y,X) + (\nabla_X Ric)(Y,Y)\right)d\tau = 2\sqrt{\bar{\tau}}\langle\nabla_Y Y, X\rangle +$$
$$\left(\int_0^{\bar{\tau}} \sqrt{\tau}\left(\nabla^2_{Y,Y}R + 2\langle R(Y,X)Y,X\rangle + 2|\nabla_Y X|^2\right)d\tau\right) +$$
$$\int_0^{\bar{\tau}} \sqrt{\tau}\left(-4(\nabla_Y Ric)(Y,X) + 2(\nabla_X Ric)(Y,Y)\right)d\tau$$

Because we know $\nabla_X X - \frac{1}{2}\nabla R + \frac{1}{2\tau}X + 2Ric(X,.) = 0$ and $Hess(f)(X,Y) = \nabla_X \nabla_Y(f) - \nabla_{\nabla_X Y}(f)$, so proof is complete ∎

**$\mathcal{L}$ and Riemannian distance**

Theorem2.5 (see[8]) let $\gamma: [0,\bar{\tau}] \to M$, $\bar{\tau} \in (0,T]$, be a $C^1$ −path starting at p and ending at q.

i) (bounding Riemannian distance by $\mathcal{L}$ ) for any $\tau \in [0,\bar{\tau}]$ we have

$$d^2_{g(0)}(p,\sigma(\tau)) \leq 2\sqrt{\tau}e^{2C_0\tau}\left(\mathcal{L}(\gamma) + \frac{2nC_0}{3}\bar{\tau}^{\frac{3}{2}}\right)$$

Where $(M^n, g(\tau)), \tau \in [0,T]$, denote a complete solution to the backward Ricci flow, and $p \in M$ shall be a base point. also we assume the curvature bound

$$\max_{(x,t)\in M\times[0,T]}\{|Rm(x,\tau)|, |Rc(x,\tau)|\} \leq C_0 < \infty$$

Remark for i) when M is noncompact. from i) we conclude for any $\bar{\tau} \in (0,T]$.

$$\lim_{q\to\infty} \bar{L}(q,\bar{\tau}) := \lim_{q\to\infty} 2\sqrt{\bar{\tau}}L(q,\bar{\tau}) = \infty$$

ii) (bounding speed at some time by $\mathcal{L}$ ) there exists $t_* \in (0,\bar{\tau})$ such that

$$\tau_*\left|\frac{\partial\gamma}{\partial\tau}(\tau_*)\right|^2_{g(\tau_*)} = \left|\frac{\partial\beta}{\partial\sigma}(\sigma_*)\right|^2_{g(\tau_*)} \leq \frac{1}{2\sqrt{\bar{\tau}}}\mathcal{L}(\gamma) + \frac{nC_0}{3}\bar{\tau}$$

Where $(\sigma) := \gamma(\tau)$, $\sigma = 2\sqrt{\tau}$, and $\sigma_* = 2\sqrt{\tau_*}$.

iii) (bounding L by Riemannian distance) for any $q \in M$ and $\bar{\tau} > 0$

$$L(q,\bar{\tau}) \leq e^{2C_0\bar{\tau}}\frac{d^2_{g(\bar{\tau})}(p,q)}{2\sqrt{\bar{\tau}}} + \frac{2nC_0}{3}\bar{\tau}^{\frac{3}{2}}$$

For proof of this theorem, we start with a remark

Remark 2.5: for $\tau_1 < \tau_2$ and $x \in M$

$$e^{-2C_0(\tau_2-\tau_1)}g(\tau_2,x) \leq g(\tau_1,x) \leq e^{2C_0(\tau_2-\tau_1)}g(\tau_2,x)$$

Proof i) suppose $\bar{\sigma} = 2\sqrt{\bar{\tau}}$ and $\beta(\bar{\sigma}) := \gamma(\bar{\tau})$. At first we compute

$$\int_0^{2\sqrt{\tau}} \left|\frac{\partial \beta}{\partial \overline{\sigma}}(\overline{\sigma})\right|^2 d\overline{\sigma} = \mathcal{L}(\gamma) - \int_{2\sqrt{\tau}}^{\sqrt{\overline{\tau}}} \left|\frac{\partial \beta}{\partial \overline{\sigma}}(\overline{\sigma})\right|^2_{g\left(\frac{\overline{\sigma}^2}{4}\right)} d\overline{\sigma} - \int_0^{\overline{\tau}} \sqrt{\overline{\tau}} R(\gamma(\overline{\tau}), \overline{\tau}) d\overline{\tau} \leq \mathcal{L}(\gamma) + \frac{2nC_0}{3}\overline{\tau}^{\frac{3}{2}}$$

Because $R \geq -nC_0$. Hence, since $g(0) \leq e^{2C_0 \tau} g(\overline{\tau})$ for $\overline{\tau} \in [0, \tau]$ we get

$$d^2_{g(0)}(p, \gamma(\tau)) \leq e^{2C_0 \tau} \left(\int_0^{2\sqrt{\tau}} \left|\frac{\partial \beta}{\partial \overline{\sigma}}(\overline{\sigma})\right|_{g\left(\frac{\overline{\sigma}^2}{4}\right)} d\overline{\sigma}\right)^2$$

$$\leq e^{2C_0 \tau} 2\sqrt{\tau} \int_0^{2\sqrt{\tau}} \left|\frac{\partial \beta}{\partial \overline{\sigma}}(\overline{\sigma})\right|^2_{g\left(\frac{\overline{\sigma}^2}{4}\right)} d\overline{\sigma} \leq 2\sqrt{\tau} e^{2C_0 \tau} \left(\mathcal{L}(\gamma) + \frac{2nC_0}{3}\overline{\tau}^{\frac{3}{2}}\right)$$

*ii) by taking $\tau = \overline{\tau}$ in proof of i) we get*

$$\frac{1}{2\sqrt{\overline{\tau}}} \int_0^{2\sqrt{\overline{\tau}}} \left|\frac{\partial \beta}{\partial \overline{\sigma}}(\overline{\sigma})\right|^2_{g\left(\frac{\overline{\sigma}^2}{4}\right)} d\overline{\sigma} \leq \frac{1}{2\sqrt{\overline{\tau}}}\mathcal{L}(\gamma) + \frac{nC_0}{3}\overline{\tau}$$

*Now with using the mean value Theorem for integrals, There exists $t_* \in (0, \overline{\tau})$ such that*

$$\left|\frac{\partial \beta}{\partial \sigma}(\sigma_*)\right|^2_{g(\tau_*)} \leq \frac{1}{2\sqrt{\overline{\tau}}}\mathcal{L}(\gamma) + \frac{nC_0}{3}\overline{\tau}$$

*iii) Let $\eta: [0, 2\sqrt{\overline{\tau}}] \to M$ be a minimal geodesic from $p$ to $q$ with respect to metric $g(\overline{\tau})$. Then because on $g(\overline{\tau})$ we have*

$$2\sqrt{\overline{\tau}} L(\gamma(\tau), \tau) = 2\sqrt{\overline{\tau}} \mathcal{L}(\gamma|_{[0,\tau]}) = \frac{\tau}{\overline{\tau}} d^2_{g(\overline{\tau})}(p, q)$$

*So*

$$L(q, \overline{\tau}) \leq \mathcal{L}(\eta) = \int_0^{2\sqrt{\overline{\tau}}} \left(\frac{\sigma^2}{4} R\left(\eta(\sigma), \frac{\sigma^2}{4}\right) + \left|\frac{d\eta}{d\sigma}\right|^2_{g(\tau)}\right) d\sigma \leq \int_0^{2\sqrt{\overline{\tau}}} \left(\frac{nC_0 \sigma^2}{4} + e^{2C_0 \overline{\tau}} \left|\frac{d\eta}{d\sigma}\right|^2_{g(\tau)}\right) d\sigma$$

$$\leq \frac{2nC_0}{3}\overline{\tau}^{\frac{3}{2}} + \frac{e^{2C_0 \overline{\tau}}}{2\sqrt{\overline{\tau}}} d^2_{g(\overline{\tau})}(p, q)$$

*So proof is complete* ■

*the $\mathcal{L}$-Jacobi equation : consider a family $\gamma(\tau, u)$ of $\mathcal{L}$-geodesics parameterized by $s$ and defined on $[\tau_1, \tau_2]$ with $0 \leq \tau_1 \leq \tau_2$. Let $Y(\tau)$ be a vector field along $\gamma$ defined by*

$$Y(\tau) = \frac{\partial}{\partial s}\gamma(\tau, s)\bigg|_{s=0}$$

*Now from second variation Theorem we obtain $Y(\tau)$ satisfies the $\mathcal{L}$ −Jacobi equation*

$$\nabla_X\nabla_X Y + R(Y,X)X - \frac{1}{2}\nabla_Y(\nabla R) + \frac{1}{2\tau}\nabla_X Y + 2(\nabla_Y Ric)(X,.) + 2Ric(\nabla_X Y,.) = 0 \quad (2.13)$$

*This is a second –order linear equation for Y .supposing that $\tau_1 > 0$ ,there is a unique vector field Y along γ solving this equation ,vanishing at $\tau_1$ with a given first-order derivative along γ at $\tau_1$.similarly ,there is a unique solution Y to this equation ,vanishing at $\tau_2$ and with a given first order derivative at $\tau_2$ .*

*Definition 2.3 a field $Y(\tau)$along an $\mathcal{L}$ −geodesic is called an $\mathcal{L}$ −jacobi field if it satisfies the $\mathcal{L}$ −jacobi equation ,equation (2.13) ,and if it vanishes at $\tau_1$.*

*Notation: for every vector field Y along γ we denote by $Jac(Y)$the expression on the left-hand side of equation(2.13) .*

*In one of remarks we considered that the vector $\lim_{\tau\to 0}\sqrt{\tau}X(\tau)$ exists .Now we give a similar result even for $\tau_1 = 0$.*

*Lemma2.5 let γ be an $\mathcal{L}-$ geodesic defined on $[0,\tau_2]$ and let $Y(\tau)$ be an $\mathcal{L}$ −Jacobi field along γ .Then $\lim_{\tau\to 0}\sqrt{\tau}\nabla_X Y$ exists ,furthermore , $Y(\tau)$ is completely determined by this limit .*

*Remark 2.6 we can consider that the bilinear pairing*

$$-\int_{\tau_1}^{\tau_2} 2\sqrt{\tau}\langle Jac(Y_1),Y_2\rangle d\tau$$

*Is a symmetric function of $Y_1$ and $Y_2$ .(here we assume that $Y_1(\tau_2) = Y_2(\tau_2) = 0$and γ is an $\mathcal{L}-$ geodesic and $Y_1,Y_2$ are vector fields along γ vanishing at $\tau_1$ ) .*

**Estimate the Hessian of the $\mathcal{L}-$ distance functions**

*Here we give an inequality for the hessian of $\mathcal{L}$ involving the integral of the vector field along.*

*Let $\gamma:[0,\bar\tau]\to M$ be an $\mathcal{L}$ −shortest curve connecting p and q .we fix a vector Y at $\tau = \bar\tau$ with $|Y|_{g_{ij}(\bar\tau)} = 1$,and extend Y along the $\mathcal{L}$ −shortest geodesic γ on $[0,\bar\tau]$ by solving the following ODE*

$$\nabla_X Y = -Ric(Y,.) + \frac{1}{2\tau}Y. \quad (2.14)$$

*Lemma2.6 (Perelman) suppose $\{Y_1,Y_2,...Y_n\}$ is an orthonormal basis at $\tau = \bar\tau$ with respect to metric $g_{ij}(\bar\tau)$ and solve for $Y(\tau)$ in the equation (2.14), then $\{Y_1(\tau),Y_2(\tau),...Y_n(\tau)\}$remains orthogonal on $[0,\bar\tau]$ and*

$$\langle Y_i(\tau),Y_j(\tau)\rangle = \frac{\tau}{\bar\tau}\delta_{ij}$$

*Proof according to (2.6) and (2.14) we get*

$$\frac{d}{d\tau}\langle Y_i,Y_j\rangle = 2Ric(Y_i,Y_j) + \langle\nabla_X Y_i,Y_j\rangle + \langle Y_i,\nabla_X Y_j\rangle$$

$$\langle Y_i(\tau),Y_j(\tau)\rangle = \frac{\tau}{\bar\tau}\delta_{ij} \Rightarrow |Y(\tau)|^2 = \frac{\tau}{\bar\tau}$$

So $\{Y_1(\tau), Y_2(\tau), \ldots Y_n(\tau)\}$ remains orthogonal on $[0, \bar{\tau}]$ with $Y_i(0) = 0$, $i = 0, 1, \ldots, n$. So proof is complete ∎

The main result of this sub section is following Theorem from Perelman.

Theorem 2.6 (Perelman[6]) suppose that $|Y|_{g_{ij}(\bar{\tau})} = 1$ at any point $q \in M$, consider an $\mathcal{L}$ −shortest geodesic $\gamma$ connecting $p$ to $q$ and extend $Y$ along $\gamma$ by solving (2.14). So the Hessian of the $\mathcal{L}$ −distance function $L$ on $M$ with $\tau = \bar{\tau}$ satisfies

$$Hess_L(Y,Y) \leq \frac{1}{\sqrt{\bar{\tau}}} - 2\sqrt{\bar{\tau}} Ric(Y,Y) - \int_0^{\bar{\tau}} \sqrt{\tau} Q(X,Y) d\tau \quad (2.15)$$

Where
$Q(X,Y) = -\nabla_Y \nabla_Y R - 2\langle R(Y,X)Y, X \rangle - 4\nabla_X Ric(Y,Y) + 4\nabla_Y Ric(Y,X) - 2Ric_\tau(Y,Y) + 2|Ric(Y,.)|^2 - \frac{1}{\tau} Ric(Y,Y)$

is the Li-Yau-Hamilton quadratic.

Proof according to the definition $Hess_L(Y,Y)$ we have
$$Hess_L(Y,Y) = Y(Y(L))(\bar{\tau}) - \langle \nabla_Y Y, \nabla L \rangle (\bar{\tau}) \quad (2.16)$$

also $Y(Y(L))(\bar{\tau}) = \delta_Y^2(L) \leq \delta_Y^2(\mathcal{L})$ .and $\nabla L(q, \bar{\tau}) = 2\sqrt{\bar{\tau}} X$. So $\langle \nabla_Y Y, \nabla L \rangle = 2\sqrt{\bar{\tau}} \langle \nabla_Y Y, X \rangle$ .Therefore from (2.16)

$$Hess_L(Y,Y) \leq \delta_Y^2(\mathcal{L}) - 2\sqrt{\bar{\tau}} \langle \nabla_Y Y, X \rangle (\bar{\tau}) \quad (2.17)$$

So from second variation formula and (2.17) we get

$$Hess_L(Y,Y) \leq \int_0^{\bar{\tau}} \sqrt{\tau} [2|\nabla_X Y|^2 + 2\langle R(Y,X)Y, X \rangle + \nabla_Y \nabla_Y R + 2\nabla_X Ric(Y,Y) - 4\nabla_Y Ric(Y,X)] d\tau$$

$$= \int_0^{\bar{\tau}} \sqrt{\tau} \left[ 2\left| -Ric(Y,.) + \frac{1}{2\tau} Y \right|^2 + 2\langle R(Y,X)Y, X \rangle + \nabla_Y \nabla_Y R + 2\nabla_X Ric(Y,Y) - 4\nabla_Y Ric(Y,X) \right] d\tau$$

So because $|Y|^2 = 1$ we get
$$\left| -Ric(Y,.) + \frac{1}{2\tau} Y \right|^2 = |Ric(Y,.)|^2 - \frac{1}{\tau} Ric(Y,Y) + \frac{1}{4\tau\bar{\tau}}$$

so
$$Hess_L(Y,Y) \leq \int_0^{\bar{\tau}} \sqrt{\tau} \left[ 2|Ric(Y,.)|^2 - \frac{2}{\tau} Ric(Y,Y) + \frac{1}{2\tau\bar{\tau}} + 2\langle R(Y,X)Y, X \rangle + \nabla_Y \nabla_Y R + 2\nabla_X Ric(Y,Y) - 4\nabla_Y Ric(Y,X) \right] d\tau$$

also because
$$\frac{d}{d\tau} Ric(Y,Y) = \frac{\partial}{\partial \tau} Ric(Y,Y) + \nabla_X Ric(Y,Y) + 2Ric(\nabla_X Y, Y)$$

and from (2.14)
$$\frac{d}{d\tau} Ric(Y,Y) = \frac{\partial}{\partial \tau} Ric(Y,Y) + \nabla_X Ric(Y,Y) - 2|Ric(Y,.)|^2 + \frac{1}{\tau} Ric(Y,Y)$$

therefore we have

$$Hess_L(Y,Y) \le \int_0^{\bar\tau} \sqrt{\tau}\left[2|Ric(Y,.)|^2 - \frac{2}{\tau}Ric(Y,Y) + \frac{1}{2\tau\bar\tau} + 2\langle R(Y,X)Y,X\rangle + \nabla_Y\nabla_Y R\right.$$
$$- 4(\nabla_Y Ric)(Y,Y)$$
$$- \left(2\frac{d}{d\tau}Ric(Y,Y) - 2\frac{\partial}{\partial\tau}Ric(Y,Y) + 4|Ric(Y,.)|^2 - \frac{2}{\tau}Ric(Y,Y)\right)$$
$$\left.+ 4\nabla_X Ric(Y,Y)\right]d\tau$$
$$= -\int_0^{\bar\tau}\left[2\sqrt{\tau}\frac{d}{d\tau}Ric(Y,Y) + \frac{1}{\sqrt{\tau}}Ric(Y,Y)\right]d\tau + \frac{1}{2\bar\tau}\int_0^{\bar\tau}\frac{1}{\sqrt{\tau}}d\tau +$$
$$+ \int_0^{\bar\tau}\sqrt{\tau}\left[2\langle R(Y,X)Y,X\rangle + \nabla_Y\nabla_Y R + \frac{1}{\tau}Ric(Y,Y) + 4\nabla_X Ric(Y,Y) - \nabla_Y Ric(X,Y)\right.$$
$$\left.+ \frac{2\partial Ric}{\partial\tau}(Y,Y) - 2|Ric(Y,.)|^2\right]d\tau = \frac{1}{\sqrt{\bar\tau}} - 2\sqrt{\bar\tau}Ric(Y,Y) - \int_0^{\bar\tau}\sqrt{\tau}Q(X,Y)d\tau$$

*So proof is complete*■

*Theorem2.7 (Perelman[8]) suppose that* $K = \int_0^{\bar\tau}\tau^{\frac{3}{2}}H(X)d\tau$ *then we have*

$$\Delta L \le \frac{n}{\sqrt{\bar\tau}} - 2\sqrt{\bar\tau}R - \frac{1}{\bar\tau}K$$

*Proof  let* $\{Y_1, Y_2, ... Y_n\}$ *be an orthonormal basis at* $\tau = \bar\tau$ *and with extend them along the shortest* $\mathcal{L}$ *−geodesic $\gamma$ and taking* $Y = Y_i$ *in (2.15)and summing over i ,we obtain*

$$\sum Hess_L(Y_i, Y_i) \le \frac{n}{\sqrt{\bar\tau}} - 2\sqrt{\bar\tau}R - \sum_{i=1}^n\int_0^{\bar\tau}\sqrt{\tau}Q(X,Y_i)d\tau$$

*But we know*   $\quad\Delta L = \sum_i Hess(L)(Y_i, Y_i)$

*So*   $\quad\Delta L \le \frac{n}{\sqrt{\bar\tau}} - 2\sqrt{\bar\tau}R - \sum_{i=1}^n\int_0^{\bar\tau}\sqrt{\tau}Q(X,Y_i)d\tau$   (2.18)

*Now we prove* $\sum_i Q(X,Y_i) = \frac{\tau}{\bar\tau}Q(X)$. *We know* $\langle Y_i(\tau), Y_i(\tau)\rangle = \frac{\tau}{\bar\tau}\delta_{ij}$ *so* $\left\{\sqrt{\frac{\bar\tau}{\tau}}Y_i(\tau)\right\}$ *is an orthonormal basis at $\tau$ .also*

$$\sum_i Q(X,Y_i) = -\frac{\tau}{\bar\tau}\Delta R + 2\frac{\tau}{\bar\tau}Ric(X,X) - 4\frac{\tau}{\bar\tau}\langle\nabla R,X\rangle + 4\frac{\tau}{\bar\tau}\sum_i\nabla_{Y_i}Ric(Y_i,X)$$
$$- 2\frac{\tau}{\bar\tau}\sum_i Ric_\tau(Y_i,Y_i) + 2\frac{\tau}{\bar\tau}\sum_i|Ric(Y_i,.)|^2 - \frac{1}{\bar\tau}\sum_i Ric(Y_i,Y_i)$$

*Tracing the second Bianchi identity gives*

$$\sum_i \nabla_{Y_i} Ric(Y_i, X) = \frac{1}{2}\langle \nabla R, X \rangle$$

*Also according to following identity*

$$Ric_\tau(Y, Y) = (\nabla R_{ij})Y^i Y^j + 2R_{ikjl}R_{kl}Y^i Y^j - 2R_{ik}R_{jk}Y^i Y^j$$ we get $\sum_i Ric_\tau(Y_i, Y_i) = \Delta R$

*Putting this together gives*

$$\sum_i Q(X, Y_i)(\tau) = \frac{\tau}{\bar{\tau}} Q(X) \quad (2.19)$$

*where* $Q(X) = -\frac{\partial R}{\partial \tau} - \frac{1}{\tau}R - 2\langle \nabla R, X \rangle + 2Ric(X, X)$

*so from (2.18),(2.19) we get*

$$\Delta L \leq \frac{n}{\sqrt{\bar{\tau}}} - 2\sqrt{\bar{\tau}}R - \int_0^{\bar{\tau}} \sqrt{\tau}\left(\frac{\tau}{\bar{\tau}}\right) Q(X) d\tau$$

$$= \frac{n}{\sqrt{\bar{\tau}}} - 2\sqrt{\bar{\tau}}R - \frac{1}{\bar{\tau}}K.$$

*Corollary 2.2 we have*

$$Hess_L(Y, Y) = \frac{1}{\sqrt{\bar{\tau}}} - 2\sqrt{\bar{\tau}}Ric(Y, Y) - \int_0^{\bar{\tau}} \sqrt{\tau}\, Q(X, Y) d\tau$$

*If and only if* $Y(\tau), \tau \in [0, \bar{\tau}]$, *is an* $\mathcal{L}$*-Jacobian field*

*Proposition 2.2(see[8]) we have*

$$\Delta L = \frac{n}{\sqrt{\bar{\tau}}} - 2\sqrt{\bar{\tau}}R - \frac{1}{\bar{\tau}}K$$

*if and only if we are on a gradient shrinking soliton with*

$$R_{ij} + \frac{1}{2\sqrt{\bar{\tau}}}\nabla_i \nabla_j L = \frac{1}{2\bar{\tau}} g_{ij}$$

*Proof when* $Y_i(\tau)\ i = 1, \ldots, n$ *are* $\mathcal{L}$*-Jacobian fields along* $\gamma$ *,we have*

$$\frac{d}{d\tau}\langle Y_i(\tau), Y_j(\tau)\rangle = 2Ric(Y_i, Y_j) + \langle \nabla_X Y_i, Y_j \rangle + \langle Y_i, \nabla_X Y_j \rangle$$

*But* $\nabla L = 2\sqrt{\bar{\tau}}X$ *so* $X = \frac{\nabla L}{2\sqrt{\bar{\tau}}}$ *therefore*

$$\frac{d}{d\tau}\langle Y_i(\tau), Y_j(\tau)\rangle = 2Ric(Y_i, Y_j) + \langle \nabla_{Y_i}\left(\frac{\nabla L}{2\sqrt{\bar{\tau}}}\right), Y_j \rangle + \langle Y_i, \nabla_{Y_j}\left(\frac{\nabla L}{2\sqrt{\bar{\tau}}}\right)\rangle$$

$$= 2Ric(Y_i, Y_j) + \frac{1}{\sqrt{\tau}} Hess_L(Y_i, Y_j)$$

but we know $\langle Y_i(\tau), Y_j(\tau) \rangle = \frac{\tau}{\bar{\tau}} \delta_{ij}$ so from last equality we get

$$R_{ij} + \frac{1}{2\sqrt{\bar{\tau}}} \nabla_i \nabla_j L = \frac{1}{2\bar{\tau}} g_{ij}$$

Also the vice versa of proof is obvious .so proof is complete ∎

**Li-Yau-Perelman distance**

*We introduce the Li-Yau-Perelman distance both on the tangent space and on space-time .The reason that the Li-Yau-Perelman distance $l = l(q, \bar{\tau})$ is easier to work with is that it is scale invariant when $\tau_1 = 0$.*

*Definition2.4 The Li-Yau-Perelman distance $l = l(q, \bar{\tau})$ is defined by*

$$l = l(q, \bar{\tau}) = \frac{L(q, \bar{\tau})}{2\sqrt{\bar{\tau}}}$$

*So in summary if we write the Perelman works on l that we proved in previous theorem, we get the following theorem.*

*Theorem2.8 (Perelman [8]) For the Li-Yau-Perelman distance $l(q, \bar{\tau})$ we have*

a) $\frac{\partial l}{\partial \bar{\tau}} = -\frac{l}{\bar{\tau}} + R + \frac{1}{2\bar{\tau}^{\frac{3}{2}}} K$

b) $|\nabla l|^2 = -R + \frac{l}{\bar{\tau}} - \frac{1}{2\bar{\tau}^{\frac{3}{2}}} K$

c) $\Delta l \leq -R + \frac{n}{2\bar{\tau}} - \frac{1}{2\bar{\tau}^{\frac{3}{2}}} K$

*Now we get the upper bound on the minimum of $l(., \tau)$ for every $\tau$*

*Lemma 2.7(see[4]) $\min_{x \in M} l(., \tau) \leq \frac{n}{2}$ for every $\tau$ .*
Proof. We know
1). $\frac{\partial L}{\partial \tau} = 2\sqrt{\tau} R - \frac{1}{2\tau} L + \frac{1}{\tau} K$
2). $\Delta L \leq \frac{n}{\sqrt{\tau}} - 2\sqrt{\tau} R - \frac{1}{\tau} K$
3). $|\nabla L|^2 = -4\tau R + \frac{2}{\sqrt{\tau}} L - \frac{4}{\sqrt{\tau}} K$
So 1,2 and 3 gives us

$$\frac{\partial l}{\partial \tau} - \Delta l + |\nabla l|^2 - R + \frac{n}{2\tau} \leq 0$$

Also 2 and 3 gives us

$$2\Delta l - |\nabla l|^2 + R + \frac{l-n}{\tau} \leq 0$$

Let $\overline{L} = 2\sqrt{\tau}L$. Therefore 1 and 2 gives us

$$\frac{\partial \overline{L}}{\partial \tau} + \Delta \overline{L} \leq 2n$$

So
$$\frac{\partial(\overline{L}-2n\tau)}{\partial \tau} + \nabla(\overline{L}-2n\tau) \leq 0.$$

Thus, by a standard maximum principle argument, $\min\{\overline{L}(q,\tau) - 2n\tau | q \in M\}$ is non-increasing and therefore $\min\{\overline{L}(q,\tau)|q \in M\} \leq 2n\tau$, so $\min l(.,\tau) \leq \frac{n}{2}$.

**Estimates on Li-Yau-Perelman distance**

*Proposition2.3 if the metrics $g_{ij}(\tau)$ have non-negative curvature operator and if the flow exists for $\tau \in [o, \tau_0]$, Then*

$$|\nabla l|^2 + R \leq \frac{cl}{\tau}$$

*For some constant c, when ever $\tau$ is bounded away from $\tau_0$, say $\tau \leq (1-c)\tau_0$, where $c > 0$.*

*Lemma2.8 (Perelman[14]) if we have a Ricci flow $\frac{\partial g_{ij}}{\partial \tau} = 2R_{ij}$, then $R(.,\tau) \geq -\frac{n}{2(\overline{\tau}-\tau)}$ whenever the flow exists for $\tau \in [0,\overline{\tau}]$.*

*Proof. We know $\frac{\partial R}{\partial \tau} = -\Delta R - 2|Ric|^2$. Look at the corresponding ODE, $\frac{\partial R}{\partial \tau} = -2|Ric|^2$. Since $R = tr(Ric)$, $|Ric|^2 \geq \frac{1}{n}R^2$ and therefore, $-2|Ric|^2 \leq \frac{-2|R|^2}{n}$, i.e. $\frac{\partial R}{\partial \tau} \leq \frac{-2|R|^2}{n}$. By solving this equation we get that the set $R(.,\tau) \geq -\frac{n}{2(\tau_0-\tau)}$ is preserved by the ODE and therefore it is preserved by the corresponding PDE.*

*Corollary2.3 For every $\overline{\tau} > 0, \exists q \in M$ and $\exists \mathcal{L}$-geodesic $\gamma: [0,\overline{\tau}] \to M$, $\gamma(0) = p$ and $\gamma(\overline{\tau}) = q$ such that $\mathcal{L}(\gamma) \leq n\sqrt{\overline{\tau}}$.*

*We proved $\min_{\overline{\tau}>0} l(.,\overline{\tau}) \leq \frac{n}{2}$. An analogy with this ideas can be found in the original proof of the Harnack inequality given by Li and Yau. They proved that, under the assumption $Ric(M) \geq -k$, a positive solution of the heat equation $\left(\frac{\partial}{\partial t} + \Delta\right)u = 0$, satisfies the gradient estimate $\frac{|gradu|^2}{u^2} - \frac{u_t}{t} \leq \frac{n}{2t}$. Along the proof of this fact, they define the function $F(x,t) = t(|grad f|^2 - f_t)$.*

*Now we give a short review Rugang Ye works on estimates for reduced lengths, that is important for the applications to have results on the Lipschitz properties on l or L.*

*Proposition2.4(see[12]) Let $\bar{\tau} \in (0,T)$. Assume that $Ric \geq -cg$ on $[0,\bar{\tau}]$ for a non-negative constant $C$. Then $L(.,\tau)$ is locally Lipschitz with respect to the metric $g(\tau)$ for each $\tau \in (0,\bar{\tau}]$. Moreover, for each compact subset $E$ of $M$, There are positive constants $A_1$ and $A_2$ such that $\sqrt{\tau}L \leq A_1$ on $E \times (0,\bar{\tau}]$ and*

$$|\dot{\gamma}(s)|^2 \leq \frac{A_2}{s}\left(1 + \frac{1}{\tau}\right)$$

*For $s \in (0,\tau]$, where $\tau \in (0,\bar{\tau}]$ and $\gamma$ denotes an arbitrary $\mathcal{L}_{0,\tau}$ −geodesic from $p$ to $q \in E$, where we denote*

$$\mathcal{L}_{a,b}(\gamma) = \int_{\sqrt{a}}^{\sqrt{b}} \left(\frac{1}{2}|\gamma'|^2 + 2Rt^2\right)dt$$

*Proposition2.5(see[1]) assume that the Ricci curvature is bounded from below on $[0,\bar{\tau}]$. Then $L(q,.)$ is locally Lipschitz on $(0,\bar{\tau}]$ for every $q \in M$. Moreover, $\tau^{\frac{3}{2}}\left|\frac{\partial L}{\partial \tau}\right|$ is bounded on $E \times (0,\bar{\tau}]$ for each compact subset $E$ of $M$.*

*Proposition2.6 assume that the Ricci curvature is bounded from below on $[0,\bar{\tau}]$. Then $L$ is locally lipschitz function on $M \times (0,\bar{\tau}]$.*

*Also similar estimates for $l$ hold in the case of bounded sectional curvature.*

*Proposition2.7(see[1]) Assume that the sectional curvature is bounded on $[0,\bar{\tau}]$. Then there is a positive constant $c = c(\tau^*)$ for every $\tau^* \in (0,\bar{\tau})$ with the following properties. For each $\tau \in (0,\tau^*]$ we have*

$$|\nabla l|^2 \leq \frac{c}{\tau}(l + \tau + 1)$$

*Almost everywhere in $M$. For each $q \in M$ we have*

$$\left|\frac{\partial l}{\partial \tau}\right| \leq \frac{c}{\tau}(l + \tau + 1)$$

*Almost everywhere in $(0,\tau^*]$.*

**Perelman's reduced Volume**

*The Perelman reduced volume is the fundamental tool which is used to establish non-collapsing results which in turn are essential in proving the existence of geometric limits. Note that reduced volume cannot be defined globally, but only on appropriate open subsets of a time-slice. But as long as one can flow an open set $U$ of a time-slice along minimizing $\mathcal{L}$ −geodesics in the direction of decreasing $\bar{\tau}$, the reduced volumes of resulting family of open sets form a monotone non-increasing function of $\bar{\tau}$. This turns out to be sufficient to extend the non-increasing results to Ricci flow with surgery. Note that Perelman's reduced*

volume resembles the expression in Huisken's monotonicity formula for the mean curvature flow. In this section we show that the reduced volume is monotonically no increasing in $\tau$ and is finite.

*Definition2.5 Given a solution of the backward Ricci flow, the reduced volume function is defined as*

$$\tilde{V}(\tau) = \int_M \tau^{\frac{-n}{2}} exp\,(-l(q,\tau))dq$$

*being $dq$ the volume form of $g(\tau)$.*

*Notice that the integrand is the heat kernel in the Euclidean space .In next propositions we summarize the main properties of this function.*

*Remark2.1(see[13]) notice that $(4\pi)^{\frac{n}{2}} = \int_M \tau^{\frac{-n}{2}} exp\,(-l(q,\tau))dq$, which is the same integral of the reduced function, but in the Euclidean space.*

*Proof. Change of variable $s = \sqrt{\tau}$ , $\frac{d\gamma}{d\tau} = \frac{1}{2\sqrt{\tau}}\frac{d\gamma}{ds}$ ,$d\tau = 2sds$ and minimize over $\gamma$ such that $\gamma(0) = p, \gamma(\bar{s}) = q$. Then we get $L(q,\bar{s}) = \frac{1}{2}\frac{d(p,q)^2}{\bar{s}}$, $L(q,\bar{\tau}) = \frac{1}{2}\frac{d(p,q)^2}{\sqrt{\bar{\tau}}}$ and so $l(q,\bar{\tau}) = \frac{1}{4}\frac{d(p,q)^2}{\bar{\tau}}$, therefore we get $\tilde{V}(\tau) = \int_{R^n} \bar{\tau}^{\frac{-n}{2}} e^{-\frac{|q|^2}{4\bar{\tau}}} dq = (4\pi)^{\frac{n}{2}}$ constant in $\bar{\tau}$.*

*In the case when $M$ is non-compact, it is not clear a priori that the integral defining the reduced volume is finite in general.*

*In fact, as the next Propositions shows, it is always finite and indeed, it is bounded above by the integral for $R^n$.*

*Theorem2.9 (Perelman [5]) the reduced volume*

$$\tilde{V}(\tau) = \int_M \tau^{\frac{-n}{2}} exp\,(-l(q,\tau))dq$$

*is non-increasing along the backward Ricci flow.*

*Proof .By applying previous Propositions we know*

$$\frac{\partial l}{\partial \tau} = R - \frac{l}{\tau} + \frac{1}{2\tau^{\frac{3}{2}}}K \quad (2.21)$$

$$|\nabla l|^2 = -R + \frac{l}{\tau} - \frac{1}{\tau^{\frac{3}{2}}}K \quad (2.22)$$

*So we get*

$$\frac{\partial l}{\partial \tau} + |\nabla l|^2 = \frac{-K}{\tau^{\frac{3}{2}}} \quad (2.23)$$

But we know
$$\Delta L|_{\bar\tau} \le \frac{n}{\sqrt{\tau}} - 2\sqrt{\tau}R - \frac{1}{\tau}K \quad (2.24)$$
So from (2.23) and (2.24) we arrive

$$\frac{\partial l}{\partial \tau} - \Delta l + |\nabla l|^2 - R + \frac{n}{2\tau} \ge 0$$

Let $\phi = \tau^{\frac{-n}{2}} \exp(-l)$. Then
$$\begin{aligned}\phi_\tau &= \left(\frac{-n}{2\tau} - \frac{\partial l}{\partial \tau}\right)\phi \\ &\le (-\Delta l + |\nabla l|^2 - R)\phi \\ &= \Delta \phi - R\phi\end{aligned}$$

This means that we have proved the following inequality

$$\phi_\tau - \Delta \phi + R\phi \le 0$$

Therefore we get

$$\frac{d}{d\tau}\tilde{V}(\tau) = \int_M (\phi_\tau + R\phi)dq \le \int_M \Delta \phi \, dq = 0$$

So $\tilde{V}(\tau)$ is non-increasing ∎

Note that equality hold if and only if $g$ is a Ricci soliton.

Definition2.6 the $\mathcal{L}-$exponential map $\mathcal{L}\exp : T_PM \times R^+ \to M$ is defined as follows: $\forall X \in T_PM$ let:
$$\mathcal{L}\exp_X(\bar\tau) = \gamma(\bar\tau)$$
Where $\gamma(\tau)$ is the $\mathcal{L}-$geodesic, starting at p and having X as the limit of $\sqrt{\tau}\dot\gamma$ ($\tau$ as $\tau \to 0$).
Denote by $J(\tau)$ the jacobian of $\mathcal{L}\exp(\tau): T_PM \to M$. Finally, we have that $\tilde{V}(\tau) = \int_{U \subset T_PM} \tau^{\frac{-n}{2}} \exp(-l(\tau))J(\tau)\, dX.$ Where $U = U_\tau = \{x \in M | \tau \le \tau(x)\}$ and $(x) = \sup\{\tau | \gamma(\tau) \text{ is a minimizing geodesic}\}$.

In next theorem we summarize the main properties of function

Theorem2.10. (see [1]):
1) If Ric is bounded from below on $[0,\tau]$ for each $\tau$, then $\tilde{V}(\tau)$ is a non-increasing function.
2) If we assume that at least one of the following conditions hold
a) Ric is bounded on $[0,\tau]$ for all $\tau$.
b) The curvature operator is non-negative for each $\tau$

Then $\tilde{V}(\tau) \le (4\pi)^{\frac{n}{2}}$ for all $\tau$.

3) If we assume that either
a) Rm is non-negative
b) The sectional curvature is bounded on $[0,\tau]$ for all $\tau$

Then we have the following equality

$$\tilde{V}(\tau_2) - \tilde{V}(\tau_1) = -\int_{\tau_1}^{\tau_2}\int_M \left(\frac{\partial l}{\partial \tau} - R + \frac{n}{2\tau}\right) e^{-l}\tau^{\frac{-n}{2}} dq d\tau.$$